\documentclass[leqno]{article}
\usepackage{amsmath,amssymb,latexsym,amscd,pst-all,mathrsfs}
\newcommand{\se}[1]{{\section{#1}} {\setcounter{equation}{0}}}
\newtheorem{theorem}{Theorem}[section]
\newtheorem{lm}{Lemma}[section]
\newtheorem{prop}{Proposition}[section]

\newtheorem{co}{Corollary}[section]

\def\k{{K\"{a}hler }}
\def\ke{{K\"{a}hler-Einstein }}

\input epsf
\begin{document}
\hbadness=10000
\title{{\bf Exponential sums, peak sections, and an alternative version of Donaldson's theorems}}
\author{Wei-Dong Ruan\\}
\maketitle

\se{Introduction}
In this paper, we provide an alternative proof (theorems \ref{da} and \ref{db}) of Donaldson's almost-holomorphic section theorem and symplectic Lefschetz pencil theorem in \cite{D1,D2}, through constructions of certain special kind of Donaldson-type sections of the line bundle based on properties of exponential sums.\\

In the original works \cite{D1,D2}, Donaldson constructed smooth symplectic hypersurfaces and Lefschetz pencils in general compact symplectic manifold $(M,\omega)$. One of the greatest strength of Donaldson's constructions lies in its applicability to ALL symplectic manifolds. These extremely important results show great promises and provide valuable tools for understanding symplectic manifolds in general. (For example, these results have been successfully applied to study symplectic 4-manifolds through subsequent works of Donaldson and many others. Due to author's ignorance of the subject, the author would like to refer the interested readers to subsequent works of Donaldson and coworkers for precise references and detail of this fascinating development.) One of the main ingredient of \cite{D1,D2} is the use of peak sections of the $U(1)$-line bundle $L$ with curvature $-i\omega$. Donaldson constructed, in a sense, the most generic almost holomorphic sections of $L^k$ that are transverse to the zero section for large $k$, which used a lot of peak sections. Such almost holomorphic sections then define the smooth symplectic hypersurfaces and Lefschetz pencils in $(M,\omega)$. (Peak section for holomorphic line bundle was also used in the work of Tian (\cite{Tian}) in \k geometry related to Yau's stability conjecture. In a paper improving Tian's result, the author made use of the peak sections for the first time in \cite{Ruan} and was fascinated by their miraculous properties.)\\

A natural question that partially motivates the current work is: How can one construct sections of $L^k$ that satisfy similar properties as Donaldson's sections using as few peak sections as possible? With fewer peak sections, the properties of individual peak sections are more visible, for such purpose, it is not only necessary to understand the behavior of a single peak section, (which is explicit and extremely nice,) it is also necessary to understand the behavior of the sum of several peak sections, (which also turns out to be surprisingly nice.) It turns out that a finite sum of peak sections can be understood in term of exponential sum, which can be viewed as generalization of polynomial functions. It will be ideal if exponential sums possess only non-degenerate isolated critical points. For such property almost implies that the pencil of hypersurfaces defined by a pencil of exponential sums will possess only isolated conic singularities. It turns out that the reality is almost as nice! When the exponents (not the coefficients) of the exponential sum are generic, the exponential sum only possesses isolated singularities with bounded multiplicity, which almost implies that the corresponding pencil of hypersurfaces will possess only isolated singularities with bounded multiplicity. With such miraculous properties of exponential sum, perturbation of sections to general position becomes a rather simple task.\\

Using such ideas from the study of exponential sums, in section 4, we are able to construct certain special kind of Donaldson-type sections without using the more involved perturbation techniques in \cite{D1,D2}. Consequently, our construction provides an alternative proof of Donaldson's almost-holomorphic section theorem and symplectic Lefschetz pencil theorem in \cite{D1,D2}.\\

We would like to illustrate the relation between peak sections and exponential sums in the simple case of $\mathbb{C}^n$ with the standard symplectic form $\omega = \sum dx_j dy_j = \frac{i}{2} \sum dz_j d \bar{z}_j$. $A = -i\alpha = \frac{1}{4} \sum (z_j d \bar{z}_j - \bar{z}_j dz_j)$ defines an $U(1)$-connection on the trivial line bundle. $-i\omega = dA$. $\sigma_0 = e^{-|z|^2/4}$ is a holomorphic peak section of the trivial line bundle that peak at 0. For a point $p \in\mathbb{C}^n$, we will use $\sigma_p$ to denote the peak section that peak at $p$.\\

Donaldson's section is essentially a linear combination of peak sections. In the $\mathbb{C}^n$ situation, let $s = \overset{l}{\underset{j=0}{\sum}} a_j \sigma_{p_j}$, where $\{p_j\}_{j=0}^l$ are points in $\mathbb{C}^n$. The important fact that we need is

\[
s = \mu(z) \sigma_0, \mbox{ where } \mu(z) = \mbox{{\small $\displaystyle\sum_{j=0}^l$}} \tilde{a}_j e^{(\bar{p}_j \cdot z)/2},\ \tilde{a}_j = a_j  e^{-|p_j|^2/4}.
\]

$\mu(z)$ is an exponential sum and $X = s^{-1}(0) = \mu^{-1}(0)$. Namely, the study of finite sum of peak sections can be reduced to the study of exponential sum.\\

Since we use much fewer peak sections than Donaldson did, our symplectic hypersurface and pencil behave quite differently. An important feature of our construction is the appearance of real \textsf{skeletons} that govern the global structure of the symplectic hypersurface and pencil. For example, Donaldson's symplectic hypersurface $X$ distribute quite evenly all over $M$. It was shown in the last section of \cite{D1} that as a current, $X$ normalized by the factor $\frac{1}{k}$ has $\omega$ as limit. In our case, to the contrary, $X$ is concentrated near a real codimension 1 skeleton $\Gamma$, (which is a stratified space we call hypercomplex with local structure similar to simplicial complex in topology. See figure 1 in page \pageref{df}.) Furthermore, as a current, the limit of $X$ normalized by the factor $\frac{1}{k}$ is a current $[\beta_\Gamma]$ that is supported on the real skeleton $\Gamma$ (theorem \ref{eg}). In a similar fashion, the pencil of symplectic hypersurfaces we construct is globally governed by the corresponding real skeleton that is a real pencil of hypercomplices in $M$ parameterized by a tree. (See figure 2 in page \pageref{de}.) The singular set of the pencil in our construction is concentrated near a real skeleton of dimension $n$. (This is the prevailing feature of our construction, where geometrically interesting sets are quite often concentrated and evenly distributed near certain real skeletons. In such way, the structure of the hypersurface and the pencil we construct can be quite explicitly understood globally through the corresponding real skeletons and locally through structure of some special exponential sums and pencils of exponential sums.)\\

One important property of Donaldson's symplectic hypersurface $X$ is that $X$ is of bounded geometry at the $k^{-\frac{1}{2}}$-scale. In our construction, there are 2 scales. The scale of bounded geometry for the real skeleton $\Gamma$ is $\epsilon$, (which is usually greater than Donaldson's scale,) the scale of bounded geometry for $X$ is $(k\epsilon)^{-1}$, (which is usually smaller than Donaldson's scale.) When $\epsilon = k^{-\frac{1}{2}}$, the 2 scales coincide and is the same as Donaldson's scale. Even in this case, our section is still somewhat less generic (using fewer peak sections) than Donaldson's section.\\

There are 2 major approaches in studying symplectic manifold $(M, \omega)$ through the associated line bundle $L$. The first is through Lefschetz pencil as initiated by Donaldson. Another is through isotropic (Lagrangian) skeleton of a section of $L$ as initiated by Eliashberg, Gromov, (\cite{EG}) and subsequently P. Biran (\cite{Biran1}) etc. A natural question that also partially motivate the current work is: What is the isotropic (Lagrangian) skeleton of a Donaldson-type section? The isotropic (Lagrangian) skeleton is more delicate to construct than the Lefschetz pencil. More understanding about the peak sections and exponential sums is needed for such purpose. We will discuss the isotropic (Lagrangian) skeleton of the section we constructed in a sequel of this work.\\

In \cite{FOR}, we devise a method to compute the Fukaya category of certain exact symplectic manifold by reducing it to the corresponding Morse category of a non-Hausdorff manifold as perturbation of the Lagrangian skeleton of the exact symplectic manifold. In a joint work \cite{BR} with A. Bondal, the method in \cite{FOR} is applied to prove Kontsevich's homological mirror symmetry conjecture for weighted projective spaces generalizing the work \cite{AKO} of Auroux, Katzarkov and Orlov for the case of weighted projective planes. It is interesting to see if the method in \cite{FOR} could be used to understand the Fukaya category of a general symplectic manifold through the Lagrangian skeleton of a Donaldson-type section.\\

{\bf Convention of notations:} $z = z^\Re + i z^\Im$ denotes the decomposition of a complex number into its real and imaginary parts. $A = O(B)$ if there is universal constant $C>0$ such that $|A| \leq C|B|$. $A\sim B$ if $A = O(B)$ and $B = O(A)$. $C$ and its variations are used as positive constants for estimates that may differ from expression to expression. $|I|$ (resp. $|\Gamma|$) denotes the cardinality of the index set $I$ (resp. the $k$-volume of the space $\Gamma$ of real dimension $k$).\\

{\bf Acknowledgement:} I would also like to thank I. Agol, H. Masur for help with Riemann surfaces and A. Libgober for help with analytic varieties.\\

\se{Generalized polynomial (exponential sum)}
Consider the map $\exp: \mathbb{C}^n \rightarrow (\mathbb{C}^*)^n$, $w = \exp(z)$. A polynomial function on $w$

\begin{equation}
\label{a}
\mu(z) = \mbox{{\small $\displaystyle\sum_{i=0}^l$}} a_i w^{m_i} = \mbox{{\small $\displaystyle\sum_{i=0}^l$}} e^{\alpha_i + m_i \cdot z}, \mbox{ where } \alpha_i = \log a_i,
\end{equation}

is also a sum of exponential functions on $z$, where the exponents are integral linear on $z$. Linear functions on $w$ correspond to exponential sums, where $l \leq n$ and $m_{[l]} :=\{m_0, \cdots, m_l\}$ is the set of vertices of a primitive $l$-simplex $\Delta$ in $\mathbb{Z}^n$.\\

We would like to consider more general exponential sum $\mu(z)$, where $m_i \in \mathbb{C}^n$, as generalization of polynomial function. For $0 \leq k \leq 2n$ and $c>0$, let $\Gamma^{(k)}$ (resp. $U_c(\Gamma^{(k)})$) contain those $z \in \mathbb{C}^n$ such that the $\mathbb{R}$-affine span of $m_{I_z} = \{m_i\}_{i\in I_z}$ (resp. $m_{I_{z,c}} = \{m_i\}_{i\in I_{z,c}}$) is of real dimension $\geq 2n-k$, where $b(z) := \max ((\alpha_0 + m_0 \cdot z)^\Re, \cdots, (\alpha_l + m_l \cdot z)^\Re)$ and $i \in I_z$ (resp. $i \in I_{z,c}$) if and only if $(\alpha_i + m_i \cdot z)^\Re$ is equal to $b(z)$ (resp. greater than $b(z)-c$). We call $\Gamma = \Gamma^{(2n-1)}$ (resp. $\Gamma^{(k)}$) the skeleton (resp. the $k$-skeleton) of the exponential sum $\mu(z)$. $U_c(\Gamma^{(k)})$ is a neighborhood of $\Gamma^{(k)}$. $\Gamma$ naturally divides $\Gamma^{(2n)} = \mathbb{C}^n$ into regions $U_j$, where $(\alpha_j + m_j \cdot z)^\Re = b(z)$. $b(z)$ can also be used to define the normalized norm $|\cdot|_p = e^{-b(p)} |\cdot|$, which is useful for doing estimates near $p \in \mathbb{C}^n$.\\

Borrowing a terminology from topology, it is reasonable to call such stratified space $\Gamma = \Gamma^{(2n-1)}$ (resp. $\Gamma^{(k)}$) a \textsf{hypercomplex} (resp. $k$-\textsf{complex}). In particular, a graph is a 1-complex.\\

$\mu(z)$ is called a \textsf{basic exponential sum} if for any $c_1>0$, there exists $c_2>0$ such that for any $z\in \Gamma^{(k)} \setminus U_{c_2} (\Gamma^{(k-1)})$ with $k\geq n$, $I_{z,c_1} = I_z$ and $\{m_i\}_{i\in I_z}$ is the vertex set of a totally real non-degenerate $(2n-k)$-simplex. Let ${\cal S}_\mu$ denote the set of such simplices. $\mu(z)$ is called \textsf{strictly basic} if the vertices of any $k$-simplex $\Delta \in {\cal S}_\mu$ with $k < n$ form a $\mathbb{C}$-linear independent set. The basic exponential sums locally are natural analogues of usual linear functions as exponential sums.\\

To quantify the (strictly) basic exponential sum conditions, we need to review some basic properties of real simplex in $\mathbb{C}^n$. Let $\Delta \subset \mathbb{C}^n$ be an $k$-simplex with vertices $\{m_0, \cdots, m_k\}$ and $k \leq n$. Assume $m_0=0$, then we may define simplex $\Delta^{\mathbb{C}}$ with vertices $\{m_0=0, m_1, Jm_1 \cdots, m_k, Jm_k\}$. Let ${\rm Vol}_{\mathbb{C}} (\Delta) := {\rm Vol} (\Delta^{\mathbb{C}})$. It is easy to see that ${\rm Vol}_{\mathbb{C}} (\Delta)$ is independent of the choice of the vertex $m_0$ of $\Delta$ necessary to define $\Delta^{\mathbb{C}}$.\\

For an $k$-simplex $\Delta \subset \mathbb{C}^n$ with $k \leq 2n$ (resp. $k \leq n$) ${\rm Diam} (\Delta) =1$, define $\delta^{\mathbb{R}} (\Delta)$ (resp. $\delta^{\mathbb{C}} (\Delta)$) to be the minimal of $[{\rm Vol} (\Delta')]^{1/l'}$ (resp. $[{\rm Vol}_{\mathbb{C}} (\Delta')]^{1/(2l')}$) for $\Delta'$ being an $l'$-face of $\Delta$ with $l' \leq k$. $\delta^{\mathbb{R}} (\Delta)$ and $\delta^{\mathbb{C}} (\Delta)$ can also be defined for general $\Delta$ by first normalizing the diameter to be 1 through scaling. $\delta^{\mathbb{C}} (\Delta) \not=0$ if and only if $\Delta$ is a totally real $k$-simplex. Namely, $\mu(z)$ is a basic exponential sum if and only if $\delta^{\mathbb{C}} (\Delta) \not=0$ for all $\Delta \in {\cal S}_\mu$.\\

For a point $m \in \mathbb{C}^n$, by requiring $\Delta'$ (in the definition of $\delta^{\mathbb{C}} (\Delta)$) to be an $l'$-face of $\Delta$ or an $l'$-simplex generated by an $(l'-1)$-face of $\Delta$ and $m$ for $l' \leq n$, we may similarly define $\delta^{\mathbb{C}}_m (\Delta)$ that is dominated by $\delta^{\mathbb{C}} (\Delta)$. For $m=0$, $\delta^{\mathbb{C}}_0 (\Delta) \not=0$ if and only if $\delta^{\mathbb{C}} (\Delta) \not=0$ and any $\{m_i\}_{i\in I}$ with $|I| \leq n$ in the vertex set of $\Delta$ is $\mathbb{C}$-linear independent. Namely, $\mu(z)$ is a strictly basic exponential sum if and only if $\delta^{\mathbb{C}}_0 (\Delta) \not=0$ for all $\Delta \in {\cal S}_\mu$.\\
\begin{lm}
\label{ae}
For any $c_1,c_2>0$, there exists $c_3>0$ such that for any simplex $\Delta \subset \mathbb{C}^n$ with ${\rm Diam} (\Delta) =1$, $\delta^{\mathbb{C}} (\Delta) \geq c_1$ and bounded distance from a point $m' \in \mathbb{C}^n$, one can find a point $m \in \mathbb{C}^n$ so that $|m -m'| \leq c_2$ and $\delta^{\mathbb{C}}_m (\Delta) \geq c_3$.
\end{lm}
{\bf Proof:} Let $\hat{B}_{c_2} (m') = \{m \in B_{c_2} (m') : \delta^{\mathbb{C}}_m (\Delta) \leq c_3 \}$. Since $\delta^{\mathbb{C}} (\Delta) \geq c_1$, it is straightforward to show that ${\rm Vol} (\hat{B}_{c_2} (m')) = O(c_3c_2^{2n-1})$. This together with ${\rm Vol} (B_{c_2} (m')) = O(c_2^{2n})$ imply that when $c_3>0$ is small enough, $B_{c_2} (m') \setminus \hat{B}_{c_2} (m')$ is non-empty.
\hfill$\Box$\\

{\bf Remark:} Lemma \ref{ae} implies that for a basic exponential sum $\mu(z)$, $e^{-m_*\cdot z} \mu(z)$ is strictly basic for some generic choice of bounded $m_* \in \mathbb{C}^n$. (Or from another perspective, $\{m_i\}_{i=0}^l \subset \mathbb{C}^n$ is shifted by $m_* \in \mathbb{C}^n$.)\\

An exponential sum $\mu(z)$ (resp. its set of exponents $m_{[l]}$) is called \textsf{strongly basic} if any $\{m_i\}_{i\in I}$ with $|I| \leq n+2$ (resp. $|I| \leq n+1$) is the vertex set of a (resp. totally) real non-degenerate $(2n-|I| +1)$-simplex. Let $\delta (m_{[l]})$ be the minimal of $\delta^{\mathbb{R}} (\Delta')$ and $\delta^{\mathbb{C}} (\Delta'')$, where $\Delta'$ (resp. $\Delta''$) is any $k$-simplex with vertices in $m_{[l]}$ and $k \leq n$ (resp. $k\leq n+1$). $\delta (m_{[l]}) \not=0$ if and only if $\mu(z)$ is strongly basic. The strongly basic condition is a condition on the set of exponents $m_{[l]}$, hence is much easier to check than the basic condition, and is very suitable for the pencil case.\\
\begin{lm}
\label{am}
Strongly basic implies basic.
\end{lm}
{\bf Proof:} Assume that $\mu(z)$ is a strongly basic exponential sum. For any $c_1>0$ and $z\in \Gamma^{(k)}$ with $k\geq n$, since the real affine span of $\{m_i\}_{i\in I_z}$ is of dimension $2n-k \leq n$, we have $\{m_i\}_{i\in I_z}$ is the vertex set of a totally real $(2n-k)$-simplex. If $I_{z,c_1} \not= I_z$, the real affine span of $\{m_i\}_{i\in I_{z,c_1}}$ is of dimension greater than $2n-k$. Hence $z \in U_{c_1} (\Gamma^{(k-1)})$.
\hfill$\Box$\\
\begin{prop}
\label{aa}
For a basic exponential sum $\mu(z)$, there exist $c, c', c''>0$ such that the zero set (resp. the critical set, resp. the critical zero set) of $\mu(z)$ is in $U_c(\Gamma)$ (resp. $U_{c'}(\Gamma^{(n)})$, resp. $U_{c''}(\Gamma^{(n-1)})$). More precisely, one can take any $c> \log l$ (resp. take $c', c''$ only depending on the geometry of $\Gamma$).
\end{prop}
{\bf Proof:} For $c> \log l$ and $z \not\in U_c(\Gamma)$, without loss of generality, we may assume $(\alpha_0 + m_0 \cdot z)^\Re = b(z)$ and $(\alpha_j + m_j \cdot z)^\Re \leq b(z) -c$ for $1 \leq j \leq l$. Then

\[
|\mu(z)| \geq |e^{\alpha_0 + m_0 \cdot z}| - \sum_{i=1}^l |e^{\alpha_i + m_i \cdot z}| \geq (e^c - l) e^{b-c} >0.
\]

Namely, the zero set of $\mu(z)$ is in $U_c(\Gamma)$.\\

For $z \not\in U_{c'}(\Gamma^{(n)})$, without loss of generality, we may assume $(\alpha_0 + m_0 \cdot z)^\Re = b(z)$ and $(\alpha_j + m_j \cdot z)^\Re \leq b(z) -c'$ for $j\not\in I$, where $|I| \leq n$. Since $\{m_i\}_{i\in I}$ are $\mathbb{C}$-linear independent, (This is the only place strictness of basic is used.) we can find $m_0^\vee$ so that $m_0 \cdot m_0^\vee =1$ and $m_j \cdot m_0^\vee =0$ for $j \in I \setminus \{0\}$. Then for $c' > c_0 = \log \left(\underset{i\not\in I}{\sum} |m_i \cdot m_0^\vee|\right)$,

\[
|m_0^\vee (\mu(z))| \geq |e^{\alpha_0 + m_0 \cdot z}| - \sum_{i\not\in I} |e^{\alpha_i + m_i \cdot z} m_i \cdot m_0^\vee| \geq (e^{c'} - e^{c_0}) e^{b-c'} >0.
\]

Namely, the critical set of $\mu(z)$ is in $U_{c'}(\Gamma^{(n)})$.\\

Assume $z \in U_{c'}(\Gamma^{(n)}) \setminus U_{c''}(\Gamma^{(n-1)})$. Without loss of generality, we may assume $(\alpha_0 + m_0 \cdot z)^\Re = b(z)$, $(\alpha_j + m_j \cdot z)^\Re \geq b(z) -c'$ for $j \leq n$ and $(\alpha_j + m_j \cdot z)^\Re \leq b(z) -c''$ for $n+1 \leq j \leq l$. Notice that

\[
e^{-m_0 \cdot z}\mu(z) = e^{\alpha_0} + \sum_{i=1}^l e^{\alpha_i + v_i \cdot z}, \mbox{ where } v_i = m_i -m_0.
\]

Since ${\rm Span}_{\mathbb{R}} (\{v_i\}_{i=1}^n)$ is totally real, $v_1, \cdots, v_n$ are complex linear independent. Let $v_0^\vee$ be the vector satisfying $v_1 \cdot v_1^\vee =1$ and $v_j \cdot v_1^\vee =0$ for $2 \leq j \leq n$. For $c'' > c'_0 = c' + \log \left(\underset{i=n+1}{\overset{l}{\sum}} |v_i \cdot v_1^\vee|\right)$,

\[
|e^{m_0 \cdot z} v_1^\vee (e^{-m_0 \cdot z}\mu(z))| \geq |e^{\alpha_1 + m_1 \cdot z}| -\!\!\! \sum_{i=n+1}^l |e^{\alpha_i + m_i \cdot z} v_i \cdot v_1^\vee| \geq e^{b-c'}(1 - e^{c'_0 -c''}) >0.
\]

Hence $z$ is not a critical point of $e^{-m_0 \cdot z}\mu(z)$. Since $\mu(z)$ and $e^{-m_0 \cdot z}\mu(z)$ have the same critical zero set, $z$ is not a critical zero point of $\mu(z)$, and the critical zero set of $\mu(z)$ is in $U_{c''}(\Gamma^{(n-1)})$.
\hfill$\Box$\\

\begin{prop}
\label{ab}
If an irreducible analytic variety $Y$ of positive dimension in $\mathbb{C}^n$ intersects $iB_R$, then it also intersects $B_R + i\partial B_R$, where $B_R$ is the closed ball of radius $R$ in $\mathbb{R}^n \subset \mathbb{C}^n$. Consequently, if $Y$ intersects $B_R + iB_R$, then it also intersects $B_{2R} + i\partial B_R$.
\end{prop}
{\bf Proof:} Through intersection with complex hyperplane, the proposition can be reduced to the situation that $Y$ is of dimension 1. If $Y$ does not intersect $B_R + i\partial B_R$, then there exists $R_1 < R$ such that

\[
Y \cap (B_R + i B_R) \subset (B_R + i B_{R_1}).
\]

Let $Y'$ be an irreducible component of $Y \cap (B_R + i\partial B_R)$ that intersects with $iB_R$. Since $\dim Y'=1$, there exists a Riemann surface $C$ with non-empty boundary and map $\phi: C \rightarrow Y'$. Let $\psi: D \rightarrow C$ be the uniformization map. Strictly speaking, $\psi$ is only defined on $D \setminus L$, where $L \subset \partial D$ is the so-called limit set. By a theorem of Alfors, when $C$ is a Riemann surface with finite genus and non-empty finitely many boundary components, which is the case in our situation, the measure of $L$ in $\partial D$ is zero.\\

Consider $f = \psi \circ \phi: D \setminus L \rightarrow \mathbb{C}^n$. We may write $f(z) = (f_1(z), \cdots, f_n(z)) = u(z) + iv(z) = (u_1(z) + iv_1(z), \cdots, u_n(z) + iv_n(z))$. By our assumption, we have $|v| \leq R_1$, $|u| \leq R$ on $D \setminus L$ and $|u| =R$ on $\partial D \setminus L$. By adjusting coordinate $z$ on $D$, we may assume $u(0)=0$. Then for any $r \in (0,1)$, we have

\[
\int_{|z|=r} (|u|^2 - |v|^2) d\theta = \bigg(\!\int_{|z|=r} f^2 \frac{dz}{iz}\bigg)^\Re\!\!\! = 2\pi (|u(0)|^2 - |v(0)|^2) = - 2\pi |v(0)|^2 \leq 0,
\]

where $f^2 = f\cdot f = \sum_{i=1}^n f_i^2$. Consequently,

\[
\int_{|z|=r} |u|^2 d\theta \leq \int_{|z|=r} |v|^2 d\theta \leq 2\pi r^2 R_1^2.
\]

Since the measure of $L$ in $\partial D$ is zero, when taking limit as $r$ approaches $1$, we have

\[
2\pi R^2 = \int_{\partial D \setminus L} |u|^2 d\theta \leq 2\pi R_1^2,
\]

which is a contradiction.
\hfill$\Box$\\
\begin{prop}
\label{ac}
A strictly basic exponential sum $\mu$ only has isolated critical points. There exists constant $N_{B_1}>0$ such that the number of critical points (counting multiplicity) in any unit ball is bounded by $N_{B_1}$. Consequently, the multiplicity of individual critical points is bounded by $N_{B_1}$. Furthermore, for any $c_1>0$, there exists $c_2>0$ such that $|d \mu (p)|_p \geq c_2$ for $p$ away from $c_1$-balls of the critical points.
\end{prop}
{\bf Proof:} The critical set is an analytic variety defined by $n$ equations. According to basic properties of analytic varieties (that can be found in \cite{GH}), for the first statement in the proposition, we only need to show that the critical set does not contain positive dimensional components.\\

Proposition \ref{aa} implies that the critical point set of $\mu$ is in a bounded neighborhood $U_c(\Gamma^{(n)})$ of the $n$-skeleton $\Gamma^{(n)}$. We only need to show that $U_c(\Gamma^{(n)})$ does not contain any positive dimensional analytic subvariety of $\mathbb{C}^n$. This is true because the $n$-skeleton $\Gamma^{(n)}$ is totally real according to our generic assumptions.\\

The actual proof will be carried out by induction. For $k \leq n$ and any $p \in \Gamma^{(k)} \setminus \Gamma^{(k-1)}$, $T_p\Gamma^{(k)}$ is a totally real $k$-dimensional subspace of $\mathbb{C}^n$. Assume the maximal distance from points in $\partial U_c(\Gamma^{(k)})$ to $\Gamma^{(k)}$ is $c_1$. It is straightforward to find a totally real $n$-dimensional vector space $V \subset \mathbb{C}^n$ containing $T_p\Gamma^{(k)}$ and $R>0$ such that ${\rm Dist} (J\partial B_R, V) >c_1$, where $B_R$ is the ball of radius $R$ in $V$. ${\cal B}_1(p) =p + (B_R + JB_R)$ and ${\cal B}_2(p) =p + (B_{2R} + JB_R)$ are neighborhoods of $p$. By taking $c' \geq c$ suitably large, one can ensure that $\partial_1 {\cal B}_2(p) = [p + (B_{2R} + J\partial B_R)]$ is outside of $U_c(\Gamma^{(k)})$ when ${\cal B}_1(p)$ is not in $U_{c'} (\Gamma^{(k-1)})$.\\

Let $Y$ be an irreducible analytic variety of positive dimension in $\mathbb{C}^n$ such that $Y \subset U_c(\Gamma^{(k)})$. We want to show that $Y$ is in effect in $U_{c'} (\Gamma^{(k-1)})$. If $Y$ is not in $U_{c'} (\Gamma^{(k-1)})$, then one can find $p \in \Gamma^{(k)} \setminus \Gamma^{(k-1)}$ such that $Y \cap {\cal B}_1(p)$ is non-empty and is not in $U_{c'} (\Gamma^{(k-1)})$. Consequently, $\partial_1 {\cal B}_2(p)$ is outside of $U_c(\Gamma^{(k)})$. By identifying $V$ isometrically with $\mathbb{R}^n \subset \mathbb{C}^n$ and $p$ with the origin of $\mathbb{C}^n$, $(p + V + JV)$ can be naturally identified with $\mathbb{C}^n$. Then proposition \ref{ab} can be applied to show that $Y \cap \partial_1 {\cal B}_2(p)$ is non-empty. This contradicts with the assumption $Y \subset U_c(\Gamma^{(k)})$. Consequently, $Y$ is in $U_{c'} (\Gamma^{(k-1)})$.\\

If $Y$ be an irreducible analytic variety of positive dimension in $\mathbb{C}^n$ such that $Y \subset U_c(\Gamma^{(n)})$, by induction, we have that $Y$ is in effect in $U_c(\Gamma^{(0)})$ for possibly bigger $c$. Since $U_c(\Gamma^{(0)})$ is a bounded set in $\mathbb{C}^n$, this is impossible and the first part of the proposition is proved.\\

For the second part of the proposition, we will first show that there exists $c>0$ such that for any point $q$, $\displaystyle \eta_{\mu,q} = \max_{r\in [1,2]} \left(\min_{z \in \partial B_r(q)} (|d\mu (z)|_q)\right) \geq c$, where recall that $|d\mu (z)|_q = e^{-b(q)} |d\mu (z)|$, $e^{b(q)} = \max ( \{|e^{\alpha_j + m_j \cdot q}|\}_{j=0}^l)$. If not so, then there exists a sequence $\{q_i\}$ such that $\lim \eta_{\mu,q_i} =0$. Let $\mu_i(z) = e^{-b(q_i)} \mu (q_i +z)$. It is easy to see that $\mu_{\infty} = \lim \mu_i$ exists and is also a non-trivial basic exponential sum. Notice that $\eta_{\mu,q_i} = \eta_{\mu_i, 0}$. $\eta_{\mu_\infty, 0} = \lim \eta_{\mu_i, 0} =0$. Consequently, $\mu_{\infty}$ has critical points on each of $\partial B_r(0)$ for $r\in [1,2]$, which contradict with the fact that $\mu_{\infty}$ has only isolated critical points.\\

Using the residue formula in \cite{GH}, the number of critical points (counting multiplicity) of $\mu$ in $D$ can be computed via the following formula:

\[
\frac{1}{C_n} \int_{\partial D} (\partial \log |d\mu|^2) \wedge (\partial \bar{\partial} \log |d\mu|^2)^{n-1}, \mbox{ where } C_n = \int_{|z|=1} \beta,
\]
\[
\beta = (\partial \log |z|^2) \wedge (\partial \bar{\partial} \log |z|^2)^{n-1} = \frac{(n-1)!}{|z|^{2n}} \sum_{i=1}^n (dz_1 d\bar{z}_1 \cdots \bar{z}_i dz_i \cdots dz_n d\bar{z}_n).
\]

Let $D = B_r(q)$ so that $r\in [1,2]$ and $\min_{z \in \partial B_r(q)} (|d\mu (z)|_q) \geq c$. Then the number of critical points (counting multiplicity) of $\mu$ in $B_r(q)$ can be computed as

\[
\frac{1}{C_n} \int_{\partial B_r(q)} (\partial \log |d\mu|^2) \wedge (\partial \bar{\partial} \log |d\mu|^2)^{n-1} \leq C \int_{\partial B_r(q)} \frac{|\nabla^2 \mu|_q^{2n-1}}{|d \mu|_q^{2n-1}},
\]

which is bounded since $|\nabla^2 \mu|_q$ is uniformly bounded on $B_r(q)$.\\

If the last part of the proposition is not true, then there exists $c_1>0$ and a sequence $\{q_i\}$ such that $B_{c_1} (q_i)$ does not contain critical point of $\mu$ and $\lim |d\mu (q_i)|_{q_i} =0$. Let $\mu_i(z) = e^{-b(q_i)} \mu (q_i +z)$. It is easy to see that $\mu_{\infty} = \lim \mu_i$ exists on $B_{c_1} (0)$ and is also a non-trivial basic exponential sum. Notice that $|d\mu_i (0)| = |d\mu (q_i)|_{q_i}$. $|d\mu_\infty (0)| = \lim |d\mu_i (0)| =0$. Consequently, $0$ is an isolated critical point of $\mu_{\infty}$, and $\mu_i$ has critical point near $0$ in $B_{c_1} (0)$ for large $i$, which is a contradiction.
\hfill$\Box$\\

Concerning the critical zero set that we will be interested, we have the following.\\
\begin{prop}
\label{ad}
For a basic exponential sum $\mu$, there exists $c_1,c_2>0$ such that $|\mu (p)|_{C^1_p} := |\mu (p)|_p + |d\mu (p)|_p \geq c_2$ outside of $U_{c_1} (\Gamma^{(n-1)})$.
\end{prop}
{\bf Proof:} Take $c_1$ to be the $2c''$ in proposition \ref{aa}. If the proposition is not true, then there exists a sequence $\{q_i\}$ outside of $U_{c_1} (\Gamma^{(n-1)})$ such that $\lim |\mu (q_i)|_{C^1_{q_i}} =0$. Let $\mu_i(z) = e^{-b(q_i)} \mu (q_i +z)$. It is easy to see that $\mu_{\infty} = \lim \mu_i$ exists and is also a non-trivial basic exponential sum. Notice that $|\mu_i (0)|_{C^1} = |\mu (q_i)|_{C^1_{q_i}}$. $|\mu_\infty (0)|_{C^1} = \lim |\mu_i (0)|_{C^1} =0$. Consequently, $0$ is an isolated critical zero point of $\mu_{\infty}$, and $\mu_i$ has critical zero point near $0$ for large $i$. Hence $\mu$ has critical zero point near $q_i$ that will be outside of $U_{c''} (\Gamma^{(n-1)})$ for large $i$, which is a contradiction to proposition \ref{aa}.
\hfill$\Box$\\

\se{Pencil of exponential sums}
Let $\mu_0$, $\mu_\infty$ be exponential sums with the same set of exponents $\{m_0, \cdots, m_l\}$.

\[
\mu_t (z) = \mu_0 + t\mu_\infty = \mbox{{\small $\displaystyle\sum_{j=0}^l$}} a_{t,j} e^{m_j \cdot z}, \mbox{ where } a_{t,j} = e^{\alpha_{0,j}} + t e^{\alpha_{\infty,j}},
\]

is a pencil of basic exponential sums, which can be guaranteed, for example, if the set of exponents satisfies the strongly basic condition. We assume $|e^{\alpha_{0,j}}| = |e^{\alpha_{\infty, j}}|$. $\mu_0$ and $\mu_\infty$ are generic if $\{e^{\alpha_{0,j} - \alpha_{\infty,j}} \}_{j=0}^l$ are well separated points in the unit circle. (For example, $\alpha_{0,j} - \alpha_{\infty,j} = \frac{2\pi j}{l+1}i$.) Let $X_t$ denote the zero set of $\mu_t$, and $Y = X_0 \cap X_\infty$ the base locus.\\

{\bf Remark:} In the pencil case, the major concern is the singular set of the pencil $\{X_t\}$ in $\mathbb{C}^n$, which is the union of singular set of $X_t$ for all $t$. The singular set of $X_t$ is the same as the critical zero set of $\mu_t$, which is unchanged if $\mu_t$ is multiplied by an exponential function. Hence, in the pencil case, we consider basic (instead of strictly basic) exponential sums.\\

Let $\alpha^\circ = (\alpha^\circ_j)_{j=0}^l \in \mathbb{R}^{l+1}$, where $\alpha^\circ_j = \alpha^\Re_{0,j} = \alpha^\Re_{\infty,j}$. For a tree $\Upsilon$ with a unique $(l+1)$-valent vertex $\tau^\circ$ and $l+1$ legs $\Upsilon_j \cong [0, -\infty)$ attach to $\tau^\circ$ for $0 \leq j \leq l$, $\{\tilde{\alpha}_\tau \}_{\tau\in \Upsilon} =$ {\footnotesize $\displaystyle \bigcup_{j=0}^l$}$\{\tilde{\alpha}_\tau \}_{\tau\in \Upsilon_j}$ is a tree in $\mathbb{R}^{l+1}$, where $\{\tilde{\alpha}_\tau \}_{\tau\in \Upsilon_j} = \alpha^\circ + [0, -\infty) e_j$ is a ray starting from $\tilde{\alpha}_{\tau^\circ} = \alpha^\circ$. Let $\tilde{\Gamma}_\tau$ be the real skeleton of the exponential sum $\overset{l}{\underset{i=0}{\sum}} e^{\tilde{\alpha}_{\tau,i} + m_i \cdot z}$ and $\tilde{\Gamma}^{(k)}$ be the union of $\tilde{\Gamma}_\tau^{(k-1)}$ for $\tau\in \Upsilon$. $\{\tilde{\Gamma}_\tau\}_{\tau \in \Upsilon}$ can be understood as a \textsf{real pencil} of real hypercomplices in $M$ (parameterized by the tree $\Upsilon$). We call the real pencil $\{\tilde{\Gamma}_\tau\}_{\tau \in \Upsilon}$ the \textsf{skeleton} of the complex pencil $\{X_t\}_{t\in \mathbb{CP}^1}$.\\
\begin{prop}
\label{ba}
The base locus of the real pencil $\{\tilde{\Gamma}_\tau\}_{\tau \in \Upsilon}$ is {\footnotesize $\displaystyle \bigcap_{\tau\in \Upsilon}$}$\tilde{\Gamma}_\tau = \Gamma^{(2n-2)}$.
\end{prop}
{\bf Proof:} When $\tau \in \Upsilon_j$ approaches infinity, all the top dimension strata of $\Gamma$ will expand in $\tilde{\Gamma}_\tau$ except those that bound $U_j$, which will move and shrink in $\tilde{\Gamma}_\tau$. (See figure 2 in page \pageref{de}.) Since every strata of $\Gamma^{(2n-2)}$ belongs to a top dimension strata of $\Gamma$ that does not bound $U_j$, we have $\Gamma^{(2n-2)} \subset \tilde{\Gamma}_\tau$. On the other hand, every top dimension strata of $\Gamma$ bound some $U_j$, and will not be in $\tilde{\Gamma}_\tau$ when $\tau \in \Upsilon_j$ approaches infinity.
\hfill$\Box$\\
\begin{prop}
\label{bb}
One can find $c>0$ such that for each $t$, there exists $\tau_t \in \Upsilon$ such that $X_t \subset U_c(\tilde{\Gamma}_{\tau_t})$.
\end{prop}
{\bf Proof:} Since $\{e^{\alpha_{0,j} - \alpha_{\infty,j}} \}_{j=0}^l$ are well separated, there exists $r_0>0$ such that $D_j = \{t\in \mathbb{C}: |t + e^{\alpha_{0,j} - \alpha_{\infty,j}}| \leq r_0\}$ are disjoint. Let $D'$ be the compliment of all such $D_j$. Define $\tau_t = \log (|t + e^{\alpha_{0,j} - \alpha_{\infty,j}}| / r_0) \in \Upsilon_j$ for $t\in D_j$ and $\tau_t = \tau^\circ$ for $t\in D'$. It is straightforward to check that $X_t \subset U_c(\tilde{\Gamma}_{\tau_t})$, where $c>0$ only depends on $l$.
\hfill$\Box$\\
\begin{prop}
\label{bc}
The singular set of the basic pencil $\{X_t\}$ is in $U_c(\tilde{\Gamma}^{(n)})$. The singular set of $Y$ is in $U_c(\Gamma^{(n-1)})$.
\end{prop}
{\bf Proof:} By proposition \ref{bb}, there exists $\tau\in \Upsilon$ such that $X_t \subset U_c(\tilde{\Gamma}_{\tau})$. Hence $\Gamma_t \subset U_{2c}(\tilde{\Gamma}_{\tau})$. Proposition \ref{aa} then implies that the singular set of $X_t$ is in $U_c(\Gamma_t^{(n-1)}) \subset U_{3c}(\tilde{\Gamma}_{\tau}^{(n-1)})$. By the definition of $\tilde{\Gamma}^{(n)}$, we have that the union of singular set of $X_t$ for all $t$ is in $U_{3c}(\tilde{\Gamma}^{(n)})$.\\

A singular point of $Y$ is always a singular point for some $X_t$. Since the singular set of $X_t$ is in $U_{3c}(\tilde{\Gamma}_{\tau}^{(n-1)})$, the singular set of $Y$ is in $U_{3c} (\Gamma^{(n-1)})$. Here we are using the fact that $\Gamma^{(2n-2)} \subset \tilde{\Gamma}_{\tau}$ implies that $\Gamma^{(2n-2)} \cap \tilde{\Gamma}_{\tau}^{(n-1)} \subset \Gamma^{(n-1)}$.
\hfill$\Box$\\
\begin{prop}
\label{bd}
A pencil of basic exponential sums $\{\mu_t\}$ only has isolated singular points. In any unit ball the number of singular points (counting multiplicity) is bounded. Consequently, the multiplicity of individual singular points is bounded. Furthermore, for any $c_1>0$, there exists $c_2>0$ such that $|\mu_t (p)|_{C^1_p} \geq c_2$ for $p$ away from $c_1$-balls of the singular points and for all $t$.
\end{prop}
{\bf Proof:} The singular set is an analytic variety. We only need to show that the singular set does not contain positive dimensional components.\\

Proposition \ref{bc} implies that the singular set of the pencil $\{X_t\}$ is in $U_c(\tilde{\Gamma}^{(n)})$. We only need to show that $U_c(\tilde{\Gamma}^{(n)})$ does not contain any positive dimensional analytic subvariety of $\mathbb{C}^n$. This can be proved in the same way as the proof of proposition \ref{ac} due to the fact that $\tilde{\Gamma}^{(n)}$ is totally real according to our generic assumptions.\\

For the second part of the proposition, we will first show that there exists $c>0$ such that for any point $q$,

\[
\eta_{\mu_t,q} = \max_{r\in [1,2]} \left(\min_{(t,z) \in \mathbb{CP}^1 \times \partial B_r(q)} \left(|\mu_t (z)|_{C^1_q}^2\right)\right) \geq c,
\]

where $|\mu_t (z)|_{C^1_q}^2 = |\mu_t (z)|_q^2 + |d\mu_t (z)|_q^2$, $|d\mu_t (z)|_q = e^{-b_t(q)} |d\mu_t (z)|$, $e^{b_t(q)} = \max ( \{|e^{\alpha_j(t) + m_j \cdot q}|\}_{j=0}^l)$. If not so, then there exists a sequence $\{q_i\}$ such that $\lim \eta_{\mu_t,q_i} =0$. Let $\mu_{i,t}(z) = \mu_{i,0}(z) + t\mu_{i,\infty}(z)$, $\mu_{i,0}(z) = e^{-b_0(q_i)} \mu_0 (q_i +z)$ and $\mu_{i,\infty}(z) = e^{-b_\infty(q_i)} \mu_\infty (q_i +z)$. It is easy to see that $\mu_{\infty,t} = \lim \mu_{i,t}$ exists and is also a non-trivial pencil of basic exponential sums. Notice that $\eta_{\mu_t,q_i} = \eta_{\mu_{i,t}, 0}$. $\eta_{\mu_{\infty,t}, 0} = \lim \eta_{\mu_{i,t}, 0} =0$. Consequently, $\mu_{\infty,t}$ has singular point on each of $\partial B_r(0)$ for $r\in [1,2]$, which contradict with the fact that $\mu_{\infty,t}$ has only isolated singular points.\\

To estimate the number of singular points, we need to use the following proposition (which is a generalization of the residue formula used in the proof of proposition \ref{ac}). This proposition should be well known as a special case of the boundary version of the well known intersection formula for divisors in term of their Chern classes. We will give a straightforward proof in line with the residue formula.\\
\begin{prop}
\label{be}
Let $f(z) = (f_1(z), \cdots, f_n(z))$ be n-tuple of holomorphic sections of a Hermitian holomorphic line bundle $(L, \|\cdot\|)$ on a complex $n$-fold $D$ with boundary $\partial D$. The number of zeros of $f$ in $D$ (assumed to be isolated) can be computed (counting multiplicity) via the following formula:

\[
\frac{1}{C_n} \left( \int_{\partial D} (\partial \log \|f\|^2) \wedge \Omega + \int_D \omega_0^n \right),
\]

where $\Omega = (\omega + \omega_0)^{n-1} + \cdots + (\omega + \omega_0) \omega_0^{n-2} + \omega_0^{n-1}$,

\[
\omega = \partial \bar{\partial} \log \|f\|^2,\ \omega_0 = - \partial \bar{\partial} \log \|f_i\|^2,\  \|f\|^2 = \|f_1\|^2 + \cdots + \|f_n\|^2.
\]
\end{prop}
{\bf Proof:} It is helpful to make $\omega_0$ vanish in a small neighborhood of the zero points first, which amounts to a trivialization of $L$ near the zero points so that $\|f_i\| = |f_i|$ near  the zero points. Since $\omega_0$ vanishes in a small neighborhood of zero points, $(\partial \log \|f\|^2) \wedge \Omega$ restricts to $(\partial \log |f|^2) \wedge (\partial \bar{\partial} \log |f|^2)^{n-1}$ near zero points. It is straightforward to check that $\omega_0^n + d((\partial \log \|f\|^2) \wedge \Omega) = (\omega + \omega_0)^n = 0$. Integrate this equation on $D$ minus small balls around zero points and use the residue formula, we get the desired result.\\

In general, we may modify $\omega_0$ to get $\omega_\epsilon$ that vanishes in a $\epsilon$-neighborhood of the zero points, such that $\omega_0 = \lim \omega_\epsilon$ away from the zero points. The formula would be true for $\omega_\epsilon$ in the place of $\omega_0$. Take limit, we get the desired formula.
\hfill$\Box$\\

We now resume the proof of proposition \ref{bd}. In our case, we take $D = \mathbb{CP}^1 \times B_r (q)$ so that $r\in [1,2]$ and $\displaystyle \min_{(t,z) \in \mathbb{CP}^1 \times \partial B_r(q)} (|\mu_t (z)|_{C^1_q}^2) \geq c$. Singular points corresponds to zero points of $(\mu_t (z), d\mu_t (z))$ that can be viewed as $(n+1)$-tuple of sections of ${\cal O}_{\mathbb{CP}^1} (1) \times {\cal O}_{B_r (q)}$. Take $\omega_0$ to be the standard Fubini-Study metric $\omega_0 = \partial \bar{\partial} \log (1 + |t|^2)$. According to proposition \ref{be}, the number of singular points (counting multiplicity) in $D$ can be computed via the following formula:

\[
\frac{1}{C_n} \int_{\partial D} (\partial \log (|\mu_t (z)|^2 + |d\mu_t (z)|^2)) \wedge \Omega \leq C \int_{\partial D} \frac{|\mu_t|_{C^2_q}^{2n-1}}{|\mu_t|_{C^1_q}^{2n-1}},
\]

which is bounded since $|\mu|_{C^2_q}$ is uniformly bounded on $D = \mathbb{CP}^1 \times B_r (q)$.\\

If the last part of the proposition is not true, then there exists $c_1>0$ and a sequence $\{q_i\}$ such that $B_{c_1} (q_i)$ does not contain singular point of the pencil and $\displaystyle \lim_{i \rightarrow +\infty} \min_{t \in \mathbb{CP}^1} (|\mu_t (q_i)|_{C^1_{q_i}}) =0$. Let $\mu_{i,t}(z) = \mu_{i,0}(z) + t\mu_{i,\infty}(z)$, $\mu_{i,0}(z) = e^{-b_0(q_i)} \mu_0 (q_i +z)$ and $\mu_{i,\infty}(z) = e^{-b_\infty(q_i)} \mu_\infty (q_i +z)$. It is easy to see that $\mu_{\infty,t} = \lim \mu_{i,t}$ exists on $B_{c_1} (0)$ and is also a non-trivial pencil of basic exponential sums. Notice that $\displaystyle \min_{t \in \mathbb{CP}^1} (|\mu_{i,t} (0)|_{C^1_0}) = \min_{t \in \mathbb{CP}^1} (|\mu_t (q_i)|_{C^1_{q_i}})$. $\displaystyle \min_{t \in \mathbb{CP}^1} (|\mu_{\infty,t} (0)|_{C^1_0}) = \lim_{i \rightarrow +\infty} \min_{t \in \mathbb{CP}^1} (|\mu_{i,t} (0)|_{C^1_0}) =0$. Consequently, $0$ is an isolated singular point of the pencil $\{\mu_{\infty,t}\}$, and the pencil $\{\mu_{i,t}\}$ has singular point near $0$ in $B_{c_1} (0)$ for large $i$, which is a contradiction.
\hfill$\Box$\\

\se{Donaldson's two theorems}
\stepcounter{subsection}
{\bf \S \thesubsection\ The symplectic hypersurface theorem:} Results in section 2 concerning exponential sums can be used to provide an alternative proof of Donaldson's symplectic hypersurface theorem in \cite{D1}. Let $(M, J, \omega, g)$ be an almost \k manifold. For suitable $R_0>0$, one can find smooth family of local complex coordinates $z_p: B_{R_0}(p) \rightarrow \mathbb{C}^n$ (up to unitary transformation) parameterized by $p \in M$ such that $z_p(p)=0$, $|J - J_p| = O(|z_p|)$ and $\omega = dz_p \wedge d\bar{z}_p$ on $B_{R_0}(p)$, where $J_p$ is the complex structure determined by $z_p$. It is straightforward to see that ${\rm Dist}_g(p_1,p_2) / |p_1 - p_2|_p  = 1 + O(\max(|p_1|, |p_2|))$ for $p_1,p_2 \in B_{R_0}(p)$, where $|p_1 - p_2|_p := |z_p(p_1) - z_p(p_2)|$. For $p,p' \in M$, the difference between $z_p$ and $z_{p'}$ (up to unitary transformation and translation) can be controlled by $O(|z_p|{\rm Dist}_g(p,p'))$.\\
\begin{prop}
\label{al}
For $\epsilon>0$, there exists a covering $\{B_\epsilon(p_i)\}_{i\in I}$ of $M$ such that ${\rm Dist}_g(p_i,p_j) \geq \epsilon$ for any $i\not= j$ and the covering is minimal in the sense that $\{B_\epsilon(p_i)\}_{i (\not= j)\in I}$ is not a covering of $M$ for any $j\in I$. $M$ can be decomposed as the union of convex polyhedron regions $\{U_i\}_{i\in I}$, where $U_i = \{p \in M: |p - p_i|_{p_i} \leq |p - p_j|_{p_j} \mbox{ for all } j\in I\}$. Furthermore $B_{\epsilon/2} (p_i) \subset U_i \subset B_\epsilon (p_i)$. (Example of $U_i$ with adjacent $U_{i_1}, \cdots, U_{i_5}$ is indicated in figure 1 in page \pageref{df}.)
\end{prop}
{\bf Proof:} We will prove by induction. Assume that we have a set $\{p_i\}_{i\in I}$ such that ${\rm Dist}_g(p_i,p_j) > \epsilon$ for any $i\not= j$ in $I$. If $\{B_\epsilon(p_i)\}_{i\in I}$ is not a covering of $M$, pick $p_{i'}$ in the compliment of $\{B_\epsilon(p_i)\}_{i\in I}$, and enlarge $I$ to include $i'$, then we still have ${\rm Dist}_g(p_i,p_j) > \epsilon$ for any $i\not= j$ in $I$. Since $\{B_{\epsilon/2} (p_i)\}_{i\in I}$ are disjoint, this induction process has to end in finite steps. We get the desired covering $\{B_\epsilon(p_i)\}_{i\in I}$ of $M$. $B_{\epsilon/2} (p_i) \subset U_i \subset B_\epsilon (p_i)$ and the minimality of the covering are easy consequences of the condition that ${\rm Dist}_g(p_i,p_j) > \epsilon$ for any $i\not= j$ in $I$.
\hfill$\Box$\\
\begin{prop}
\label{af}
There exists $c_2>0$ such that when $p_I = \{p_i\}_{i\in I}$ is chosen generically, $\delta (\Delta) \geq c_2$ for any $\Delta \in \tilde{{\cal S}}_\epsilon (I)$, where $\tilde{{\cal S}}_\epsilon (I)$ denotes the set of simplex with vertices in $\{p_i\}_{i\in I}$ and lengths less or equal to $2\epsilon$ 1-edges.
\end{prop}
{\bf Proof:} Fix $c_1>0$ small. We will prove by inductive construction of $I_1, \cdots, I_N$. $\{p_i\}_{i\in I_1}$ can be constructed according to proposition \ref{al} so that ${\rm Dist}_g(p_i,p_j) \geq 4\epsilon$ for any $i\not= j$ in $I_1$ and $\{B_{4\epsilon}(p_i)\}_{i\in I_1}$ forms a covering of $M$. $I = I_1$ trivially satisfies ${\rm Dist}_g(p_i,p_j) > (1-c_1)\epsilon$ for any $i\not= j$ in $I$, $\delta (\Delta) \geq 2c_2$ for any $\Delta \in \tilde{{\cal S}}_\epsilon (I)$. Then for any $l$-simplex $\Delta \in \tilde{{\cal S}}_\epsilon (I)$, ${\rm Vol} (\Delta) \geq (c_2 \epsilon)^l$ (resp. ${\rm Vol}_{\mathbb{C}} (\Delta) \geq (c_2 \epsilon)^{2l}$) when $l\leq n+1$ (resp. $l\leq n$).\\

Let $\langle \Delta \rangle$ (resp. $\langle \Delta \rangle_{\mathbb{C}}$) be the minimal real (resp. complex) affine space containing $\Delta$. Let ${\cal S}_\epsilon^{\mathbb{C}} (I)$ (resp. ${\cal S}_\epsilon (I)$) denote the set of $l$-simplex $\Delta \in \tilde{{\cal S}}_\epsilon (I)$ with $l\leq n-1$ (resp. $l\leq n$). For $\Delta \in {\cal S}_\epsilon^{\mathbb{C}} (I)$ (resp. $\Delta \in \tilde{{\cal S}}_\epsilon (I)$) we may define

\[
\hat{B}^{\mathbb{C}}_{2\epsilon} (\Delta) = B_{2\epsilon} (\Delta) \cap B_{c_3\epsilon} (\langle \Delta \rangle_{\mathbb{C}}) ({\rm resp.}\  \hat{B}_{2\epsilon} (\Delta) = B_{2\epsilon} (\Delta) \cap B_{c_3\epsilon} (\langle \Delta \rangle)).
\]

Clearly, $|\hat{B}^{\mathbb{C}}_{2\epsilon} (\Delta)| \leq Cc_3 |B_{2\epsilon} (\Delta)|$ (resp. $|\hat{B}_{2\epsilon} (\Delta)| \leq Cc_3 |B_{2\epsilon} (\Delta)|$).\\

If $\{B_\epsilon(p_i)\}_{i\in I}$ is not a covering of $M$, by taking $c_3>0$ small (depending $c_1$), one can find $p_{i'}$ in the compliment of $\{B_{(1-c_1) \epsilon} (p_i)\}_{i\in I}$, $\{\hat{B}^{\mathbb{C}}_{2\epsilon} (\Delta)\}_{\Delta \in {\cal S}_\epsilon^{\mathbb{C}} (I)}$ and $\{\hat{B}_{2\epsilon} (\Delta)\}_{\Delta \in {\cal S}_\epsilon (I)}$. More precisely, first find $p'_{i'}$ in the compliment of $\{B_{\epsilon} (p_i)\}_{i\in I}$. Then $B_{c_1 \epsilon} (p'_{i'})$ is in the compliment of $\{B_{(1-c_1) \epsilon} (p_i)\}_{i\in I}$. It is easy to see that the number of $\Delta \in {\cal S}_\epsilon^{\mathbb{C}} (I)$ (resp. $\Delta \in \tilde{{\cal S}}_\epsilon (I)$) such that

\[
\hat{B}^{\mathbb{C}}_{2\epsilon} (\Delta) \cap B_{c_1 \epsilon} (p'_{i'}) \not= \emptyset\ ({\rm resp.}\  \hat{B}_{2\epsilon} (\Delta) \cap B_{c_1 \epsilon} (p'_{i'}) \not= \emptyset)
\]

is bounded. Then

\[
\left|B_{c_1 \epsilon} (p'_{i'}) \setminus \left[ \! \left( \bigcup_{\Delta \in {\cal S}_\epsilon^{\mathbb{C}} (I)} \!\!\! \hat{B}^{\mathbb{C}}_{2\epsilon} (\Delta)\right) \cap \left( \bigcup_{\Delta \in \tilde{{\cal S}}_\epsilon (I)} \!\!\! \hat{B}_{2\epsilon} (\Delta) \right)\!\right]\right| \geq C(c_1\epsilon)^{2n} -O(c_3\epsilon^{2n}) >0
\]

when $c_3$ is small. Hence we can find the desired $p_{i'}$. Let $(\Delta, p_{i'})$ denote the simplex generated by $\Delta$ and $p_{i'}$. For $l$-simplex $\Delta \in {\cal S}_\epsilon^{\mathbb{C}} (I)$ (resp. $\Delta \in {\cal S}_\epsilon (I)$), our choice of $p_{i'}$ ensure that ${\rm Vol} (\Delta, p_{i'}) \geq c_3 \epsilon(c_2 \epsilon)^l/(2n)$ (resp. ${\rm Vol}_{\mathbb{C}} (\Delta, p_{i'}) \geq (c_3 \epsilon)^2 (c_2 \epsilon)^{2l}/(2n)^2$) when $l\leq n+1$ (resp. $l\leq n$).\\

Let $I_2$ be the set of such $i'$ so that ${\rm Dist}_g(p_i,p_j) \geq 4\epsilon$ for any $i\not= j$ in $I_2$ and $I_2$ is maximal. Enlarge $I$ to include $I_2$, we still have ${\rm Dist}_g(p_i,p_j) > (1-c_1) \epsilon$ for any $i\not= j$ in $I$. By taking $c_2$ small (depending on $c_3$), we still have ${\rm Vol} (\Delta) \geq (c_2 \epsilon)^l$ (resp. ${\rm Vol}_{\mathbb{C}} (\Delta) \geq (c_2 \epsilon)^{2l}$) when $l\leq n+1$ (resp. $l\leq n$) for any $l$-simplex $\Delta \in \tilde{{\cal S}}_\epsilon (I)$ with the enlarged $I$. Consequently, $\delta (\Delta) \geq 2c_2$ for any $\Delta \in \tilde{{\cal S}}_\epsilon (I)$.\\

Through induction, we may construct $I = I_1 \cup \cdots \cup I_N$. Since $\{B_{(1-c_1) \epsilon/2} (p_i)\}_{i\in I}$ are disjoint, this induction process has to end in finite steps. We get the desired covering $\{B_\epsilon(p_i)\}_{i\in I}$ of $M$.
\hfill$\Box$\\

Suppose that there is a line bundle $L$ with an $U(1)$-connection that has curvature $-i\omega$. For suitable trivialization of $L^k$ on $B_{R_0}(p)$, the connection 1-form $A = \frac{k}{4} (z_p d \bar{z}_p - \bar{z}_p dz_p)$ satisfying $-ik\omega = dA$. Under such trivialization, $\sigma_p = \rho_p e^{-k|z_p|^2/4}$ represents a peak section of the line bundle $L^k$ that peaks at $p$ and is supported in $B_{R_0}(p)$, where $\rho_p$ is a smooth cut-off function that equals to 1 on $B_{R_0/2}(p)$ and equals to 0 outside of $B_{R_0}(p)$. Let $X = s^{-1}(0)$, where $s = \underset{j\in I}{\sum} a_j \sigma_{p_j}$, $|a_j|=1$.\\

{\bf Remark:} The same arguments will also apply when $\log |a_j|$ is uniformly bounded, where we need to modify the definition of $U_i$ and $\Gamma$ slightly. More precisely, $U_i = \{p \in M: \frac{4}{k} (\log |a_i| - \log |a_j|) + |p - p_i|_{p_i}^2 \leq |p - p_j|_{p_j}^2 \mbox{ for all } j\in I\}$ in general. The assumption $|a_j|=1$ makes the notation much simpler.\\

For $l \geq n-1$, let $\Gamma^{(l)} \subset M$ contains those $p \in M$ such that more than $2n-l$ elements in $\{\log |\sigma_{p_j} (p)|\}_{j\in I}$ is equal to $b(p)$, where $b(z) := \max (\{\log |\sigma_{p_j} (z)|\}_{j\in I})$. We call $\Gamma = \Gamma^{(2n-1)}$ (resp. $\Gamma^{(l)}$) the skeleton (resp. the $l$-skeleton) of $X$ or $s$. Notice that $\Gamma^{(2n)} = M$ and $\Gamma =$ {\footnotesize $\displaystyle \bigcup_{i\in I}$}$\partial U_i$. (See figure 1.) $b(z)$ can also be used to define the normalized norm $|\cdot|_p = e^{-b(p)} |\cdot|$, which is useful for doing estimates near $p \in M$.\\

\pspicture(-6,-2.5)(6,1.7)
\rput(-3,0){
\psset{linecolor=lightgray,linestyle=dashed}
\gray
\rput(-.3,.4){$U_i$}
\rput(-.1,1.3){$U_{i_1}$}
\rput(-1.3,0){$U_{i_2}$}
\rput(-.6,-1.2){$U_{i_3}$}
\rput(.7,-1.2){$U_{i_4}$}
\rput(1.3,.5){$U_{i_5}$}
\psline(.7,1)(1.3,1.65)
\psline(-.74,.95)(-1.3,1.6)
\psline(1.02,-.4)(1.85,-.65)
\psline(-1,-.5)(-1.8,-.7)
\psline(0,-1.2)(0,-1.7)
\rput(-1.9,-.7){$\Gamma$}
\pspolygon(.7,1)(-.74,.95)(-1,-.5)(0,-1.2)(1.02,-.4)

\psset{linecolor=black,linestyle=solid}
\black
\pscircle[linewidth=.5pt](0,0){1.5}
\psarc[linewidth=.5pt](0,2){1.5}{200}{340}
\psarc[linewidth=.5pt](1.8,.4){1.5}{110}{270}
\psarc[linewidth=.5pt](-1.8,.3){1.5}{270}{60}
\psarc[linewidth=.5pt](-1.3,-1.5){1.5}{350}{110}
\psarc[linewidth=.5pt](1.3,-1.5){1.5}{60}{190}
\qdisk(0,0){1.5pt}
\rput(.2,-.2){$p_i$}
\qdisk(0,2){1.5pt}
\rput(.3,1.8){$p_{i_1}$}
\qdisk(-1.8,.3){1.5pt}
\rput(-1.6,.6){$p_{i_2}$}
\qdisk(-1.3,-1.5){1.5pt}
\rput(-1.3,-1.2){$p_{i_3}$}
\qdisk(1.3,-1.5){1.5pt}
\rput(1.6,-1.3){$p_{i_4}$}
\qdisk(1.8,.4){1.5pt}
\rput(1.9,.1){$p_{i_5}$}}

\rput(3,0){
\pscircle[linewidth=.5pt,linestyle=dashed,linecolor=lightgray](0,0){1.5}
\psarc[linewidth=.5pt,linestyle=dashed,linecolor=lightgray](0,2){1.5}{200}{340}
\psarc[linewidth=.5pt,linestyle=dashed,linecolor=lightgray](1.8,.4){1.5}{110}{270}
\psarc[linewidth=.5pt,linestyle=dashed,linecolor=lightgray](-1.8,.3){1.5}{270}{60}
\psarc[linewidth=.5pt,linestyle=dashed,linecolor=lightgray](-1.3,-1.5){1.5}{350}{110}
\psarc[linewidth=.5pt,linestyle=dashed,linecolor=lightgray](1.3,-1.5){1.5}{60}{190}
\psline(.7,1)(1.3,1.65)
\psline(-.74,.95)(-1.3,1.6)
\psline(1.02,-.4)(1.85,-.65)
\psline(-1,-.5)(-1.8,-.7)
\psline(0,-1.2)(0,-1.7)
\qdisk(0,0){1.5pt}
\rput(.2,-.2){$p_i$}
\rput(-.3,.4){$U_i$}

\psset{linecolor=gray}
\gray
\qdisk(0,2){1pt}
\rput(.3,1.8){$p_{i_1}$}
\rput(-.1,1.3){$U_{i_1}$}
\qdisk(-1.8,.3){1pt}
\rput(-1.6,.6){$p_{i_2}$}
\rput(-1.3,0){$U_{i_2}$}
\qdisk(-1.3,-1.5){1pt}
\rput(-1.3,-1.2){$p_{i_3}$}
\rput(-.6,-1.2){$U_{i_3}$}
\qdisk(1.3,-1.5){1pt}
\rput(1.6,-1.3){$p_{i_4}$}
\rput(.7,-1.2){$U_{i_4}$}
\qdisk(1.8,.4){1pt}
\rput(1.9,.1){$p_{i_5}$}
\rput(1.3,.5){$U_{i_5}$}
\psset{linecolor=black}
\black
\rput(-1.9,-.7){$\Gamma$}
\pspolygon(.7,1)(-.74,.95)(-1,-.5)(0,-1.2)(1.02,-.4)}
\stepcounter{figure}
\label{df}
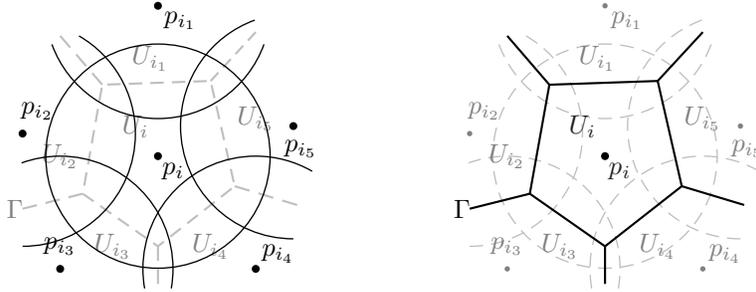
\uput{2}[d](0,0){Figure \thefigure: The ball covering and the real skeleton $\Gamma$}
\endpspicture

\begin{prop}
\label{ak}
When $\epsilon$ is small and $\{p_i\}_{i\in I}$ satisfies the generic condition in proposition \ref{af}, for $2n \geq l \geq n$, $\Gamma^{(l)}\setminus \Gamma^{(l-1)}$ is smooth of dimension $l$ with tangent space at each point containing only complex subspace of dimension less than or equal to $l-n$. Locally, each component of $\Gamma^{(l)}\setminus \Gamma^{(l-1)}$ is the common boundary of $(2n-l+1)$ components of $\Gamma^{(l+1)}\setminus \Gamma^{(l)}$ for $2n-1 \geq l \geq n$.
\hfill$\Box$
\end{prop}

As in the exponential sum case, it is reasonable to call such stratified space $\Gamma = \Gamma^{(2n-1)}$ (resp. $\Gamma^{(l)}$) a \textsf{hypercomplex} (resp. $l$-\textsf{complex}) in $M$.\\

{\bf Remark:} $\Gamma^{(l)}$ for $l<n-1$ that is unnecessary for our arguments, can be similarly defined, with structure much more complicated and less useful under our generic condition on $\{p_i\}_{i\in I}$.\\

For $p \in X$, we need to use the coordinate $z = \epsilon k z_p$. (From now on, the distance would respect the metric $\tilde{g} = \epsilon^2 k^2 g$. We will use $\tilde{B}_r$ to denote the radius $r$ ball with respect to $\tilde{g}$. For example, $\tilde{B}_{r}(p)$ would amount to $B_{\frac{r}{k\epsilon}}(p)$ under $g$.) Under the coordinate $z$, $A = -ik\alpha = \frac{1}{4\epsilon^2 k} (z d \bar{z} - \bar{z} dz)$, $-ik\omega = dA = \frac{1}{2\epsilon^2 k} dz d \bar{z}$, $\epsilon^2 k^2 \omega = dx dy$ and $\sigma_p = \rho_p e^{-|z|^2/(4\epsilon^2 k)}$. Define $\mu_p$ by $s = \mu_p \sigma_p$, we have\\
\begin{prop}
\label{ag}
For $p \in X$ with suitable choice of $m_*$, we may write $\mu_p = e^{m_* \cdot z} \mu^\circ_p + \mu'_p$ so that

\[
\mu^\circ_p = \mbox{{\small $\displaystyle\sum_{i\in I_p}$}} a_i e^{-\epsilon^2 k|m_i|^2} e^{(m_i-m_*) \cdot z}
\]

is a strictly basic exponential sum, where $I_p := \{i\in I: {\rm Dist}_g(p,p_i) \leq 2\epsilon\}$. Furthermore, $|\mu'_p|_{C^1_p} = O(c_{\epsilon,k})$ when $|z|$ is bounded, where $c_{\epsilon,k} = \max(\epsilon, \frac{1}{\epsilon k}, e^{-c\epsilon^2k})$ for certain fixed $c>0$.
\end{prop}
{\bf Proof:} Proposition \ref{af} implies that $\mu^\circ_p$ is strongly basic. Then with suitable choice of $m_*$, $\mu^\circ_p$ is strictly basic. The only thing remains to be shown is that $|\mu'_p|_{C^1_p} = O(c_{\epsilon,k})$ when $|z|$ is bounded.\\

Let $z_i = \epsilon k (z_{p_i} - z_{p_i} (p))$ then $z_i(p)=0$. Since ${\rm Dist}_g (p,p_i)=O(\epsilon)$, modify $z_i$ suitably by unitary transformation, we may further assume $|z_i - z| = O(\epsilon |z|)$. Let $A_i$ and $\sigma_{p,i}$ be the connection 1-form and peak section defined using $z_i$. Then $\sigma_{p,i} = e^{u_i + iv_i} \sigma_p$, where $u_i = (|z_i|^2 -|z|^2)/(4\epsilon^2 k)$,  $A_i - A = idv_i$. Notice that $u_i, v_i = O(|z|^2/ (\epsilon k))$.

\[
\sigma_{p_i} = e^{-(\epsilon^2 k|a_i|^2 - 2(\bar{a}_i \cdot z_i))/4} \sigma_{p,i} = e^{-\epsilon^2 k|m_i|^2} e^{m_i \cdot z} e^{w_i} \sigma_p.
\]

where $z_i(p_i) = \epsilon^2 k a_i$, $m_i = \bar{a}_i/2$ and $w_i = m_i \cdot (z_i -z) + u_i + iv_i$.\\

It is straightforward to see that $|w_i|_{C^1} = O(\epsilon, 1/(\epsilon k))$ and $|\sigma_{p_j}/\sigma_p|_{C^1_p} = O(e^{-c\epsilon^2k})$ when $|z|$ is bounded and $j\not\in I_p$. Hence $|\mu'_p|_{C^1} = O(c_{\epsilon,k})$.
\hfill$\Box$\\
\begin{prop}
\label{ai}
For any $C_1>0$ there exists $0< C_2 < C_1$ such that for $c_{\epsilon,k}$ small, $|\mu_{p'_1}|_{C^1_{p'_2}} (p) \leq C_2$ implies that  $|\mu_{p''_1}|_{C^1_{p''_2}} (p) \leq C_1$ and vis versa, for $p,p'_1,p'_2,p''_1,p''_2 \in M$ with bounded mutual distances under $\tilde{g}$.
\end{prop}
{\bf Proof:} Since $|\mu_{p'_1}|_{C^1_{p'_2}} (p) = e^{-b(p'_2)} |\mu_{p'_1}|_{C^1} (p)$, $\mu_{p'_1} (p) = e^{b'(p)} \mu_{p''_1} (p)$, where $b'(p) = [\log (\sigma_{p''_1}/ \sigma_{p'_1})] (p)$, the proposition is a consequence of the boundedness of $|b(p'_2) - b(p''_2)|$ and $|b'(p)|_{C^1}$ under the assumption of the proposition.
\hfill$\Box$\\
\begin{prop}
\label{ah}
For $c_{\epsilon,k}$ small, there exist $R_1,C_3>0$, a set of points $\{q'_i\}_{i\in \tilde{I}} \subset X$ and $\{\gamma^\circ_{q'_i}\}_{i\in \tilde{I}}$, where $\gamma^\circ_{q'_i}$ is a set of critical points of $\mu^\circ_{q'_i}$, so that $\{\tilde{B}_{3R_1} (\gamma^\circ_{q'_i})\}_{i\in \tilde{I}}$ are disjoint and $|d\mu_p (p)|_p \geq C_3$ for $p \in X \setminus \tilde{B}_{R_1} (\gamma^\circ_{\tilde{I}})$, where $\gamma^\circ_{\tilde{I}} = \bigcup_{i \in \tilde{I}} \gamma^\circ_{q'_i}$.
\end{prop}
{\bf Proof:} Let $R_1$ satisfy $16R_1 N_{B_1} <1$. Proposition \ref{ac} implies that there exists $C_1>0$ such that away from $R_1$-balls centered at critical points of $\mu^\circ_p$, $|d\mu^\circ_p|_p \geq C_1$. Proposition \ref{ag} implies that for $c_{\epsilon,k}$ small, there exists a constant $C_2>0$ such that $|d\mu_p (z)|_p \leq C_2$ implies that $|d\mu^\circ_p (z)|_p \leq C_1$.\\

Let $q'_0 \in X$ be a point satisfying $|d\mu_{q'_0} (q'_0)|_{q'_0} \leq C_3 < C_2$, then $|d\mu^\circ_{q'_0} (q'_0)| \leq C_1$. By proposition \ref{ac}, $q'_0$ is in the $R_1$-ball of a critical point $q_0$ of $\mu^\circ_{q'_0}$. Let $\gamma_{q'_0}$ be the maximal connected graph in $\tilde{B}_1(q'_0)$ with the vertex set $\gamma^\circ_{q'_0} = \{q_{0,j}\}_{j\in A_0}$ (including $q_{0,0} = q_0$) containing critical points of $\mu^\circ_{q'_0}$ such that any leg of $\gamma_{q'_0}$ has length less than $8R_1$ and $\displaystyle \min_{q\in \tilde{B}_{R_1} (q_{0,j})} (|\mu_{q}|_{C^1_q} (q)) \leq C_3$. The condition $16R_1 N_{B_1} <1$ implies that $\gamma_{q'_0} \subset \tilde{B}_{1/2} (q'_0)$.\\

Take $q'_1 \in X \setminus \tilde{B}_{R_1} (\gamma^\circ_{q'_0})$ satisfying $|d\mu_{q'_1} (q'_1)|_{q'_1} \leq C_3 < C_2$, by similar procedure we may construct $\gamma_{q_1}$ with vertices $\{q_{1,j}\}_{j\in A_1}$ (including $q_{1,0} = q_1 \in \tilde{B}_{R_1} (q'_1)$) that satisfy the additional condition: $q_{1,j} \not\in \tilde{B}_{2R_1} (\gamma^\circ_{q'_0})$ for $j\in A_1$. For $j_1\in A_1$, by definition, there exists $q \in \tilde{B}_{R_1} (q_{1,j_1})$ such that $|\mu_{q}|_{C_q^1} (q) \leq C_3$. Then by proposition \ref{ai}, when $C_3$ is small, $|\mu_{q'_0}|_{C_{q'_0}^1} (q) \leq C_2$, $|d\mu^\circ_{q'_0} (q)|_{q'_0} \leq C_1$. By proposition \ref{ac}, $q$ is in the $R_1$-ball of a critical point $q_{0,j_0}$ of $\mu^\circ_{q'_0}$ and $q_{1,j_1} \in \tilde{B}_{2R_1} (q_{0,j_0})$, which together with $q_{1,j_1} \not\in \tilde{B}_{2R_1} (\gamma^\circ_{q'_0})$ implies that $j_0 \not\in A_0$. Hence $q_{0,j_0} \not\in \tilde{B}_{8R_1} (\gamma^\circ_{q'_0})$, $q \not\in \tilde{B}_{7R_1} (\gamma^\circ_{q'_0})$ and $q_{1,j_1} \not\in \tilde{B}_{6R_1} (\gamma^\circ_{q'_0})$. Consequently, $\tilde{B}_{3R_1} (\gamma^\circ_{q'_0})$ is disjoint from $\tilde{B}_{3R_1} (\gamma^\circ_{q'_1})$.\\

It remains to verify that $q_{1,0} = q_1 \not\in \tilde{B}_{2R_1} (\gamma^\circ_{q'_0})$. Since $|\mu_{q'_1}|_{C_{q'_1}^1} (q'_1) = |d\mu_{q'_1} (q'_1)|_{q'_1} \leq C_3$, by proposition \ref{ai}, when $C_3$ is small, $|\mu_{q'_0}|_{C_{q'_0}^1} (q'_1) \leq C_2$, $|d\mu^\circ_{q'_0} (q'_1)|_{q'_0} \leq C_1$. By proposition \ref{ac}, $q'_1$ is in the $R_1$-ball of a critical point $q_{0,j_0}$ of $\mu^\circ_{q'_0}$, which together with $q'_1 \not\in \tilde{B}_{R_1} (\gamma^\circ_{q'_0})$ imply that $j_0 \not\in A_0$. Hence $q_{0,j_0} \not\in \tilde{B}_{8R_1} (\gamma^\circ_{q'_0})$, which together with $q_{1,0} \in \tilde{B}_{2R_1} (q_{0,j_0})$ imply that $q_{1,0} \not\in \tilde{B}_{2R_1} (\gamma^\circ_{q'_0})$.\\

Through induction, we can construct $\{\gamma_{q'_i}\}_{i\in \tilde{I}}$ so that $\{\tilde{B}_{3R_1} (\gamma^\circ_{q'_i})\}_{i\in \tilde{I}}$ are disjoint and $|d\mu_p (p)| \geq C_3$ for $p \in X \setminus \tilde{B}_{R_1} (\gamma^\circ_{\tilde{I}})$, where $\gamma^\circ_{\tilde{I}} = \bigcup_{i \in \tilde{I}} \gamma^\circ_{q'_i}$.
\hfill$\Box$\\

For each $i\in \tilde{I}$, we can construct a cutoff function $\tilde{\rho}_i$ (resp. $\hat{\rho}_i$) such that $\tilde{\rho}_i =1$ (resp. $\hat{\rho}_i =1$) on $\tilde{B}_{2R_1} (\gamma^\circ_{q'_i})$ (resp. $\tilde{B}_{R_1} (\gamma^\circ_{q'_i})$) and $\tilde{\rho}_i =0$ (resp. $\hat{\rho}_i =0$) away from $\tilde{B}_{3R_1} (\gamma^\circ_{q'_i})$ (resp. $\tilde{B}_{2R_1} (\gamma^\circ_{q'_i})$). Let $\hat{X} = \hat{s}^{-1}(0)$ and $\tilde{X} = \tilde{s}^{-1}(0)$, where

\[
\tilde{s} = s - \mbox{{\small $\displaystyle \sum_{i\in \tilde{I}}$ }} \tilde{\rho}_i \mu'_{q'_i} \sigma_{q'_i}, \ \hat{s} = \tilde{s} + \mbox{{\small $\displaystyle \sum_{i\in \tilde{I}}$ }} \hat{\epsilon}_i \hat{\rho}_i e^{m_* \cdot z} \sigma_{q'_i},
\]

$\hat{\epsilon}_i$ is a suitable small constant such that $\mu^\circ_{q'_i} + \hat{\epsilon}_i$ has no critical zero points in $\tilde{B}_{R_1} (\gamma_{q_i})$.\\
\begin{theorem}
\label{dc}
When $c_{\epsilon,k}$ is small enough, there exists $c>0$ such that $X$ is in $B_{\frac{c}{k\epsilon}} (\Gamma)$, and $X$ is smooth of bounded geometry at the $\frac{1}{k \epsilon}$-scale outside of $B_{\frac{c}{k\epsilon}} (\Gamma^{(n-1)})$.
\end{theorem}
{\bf Proof:} Near $p\in X$, under the coordinate $z=\epsilon kz_p$ (or the metric $\tilde{g} = \epsilon^2 k^2 g$), by proposition \ref{ag}, $\Gamma$ (resp. $\mu_p$, resp. $X$) is an $O(c_{\epsilon,k})$-perturbation of $\Gamma_p$ (resp. $\mu^\circ_p$, resp. $X_p$), where $\Gamma_p$ is the skeleton of $X_p = (\mu^\circ_p)^{-1} (0)$. When $c_{\epsilon,k}$ is small enough, $B_{\frac{c}{k\epsilon}} (\Gamma^{(l)}) = \tilde{B}_c (\Gamma^{(l)}) \sim U_c (\Gamma_p^{(l)})$ near $p$ for $l\geq n-1$. The theorem is then a consequence of proposition \ref{aa} and (proposition \ref{ad} for bounded geometry part).
\hfill$\Box$\\
\begin{theorem}
\label{dd}
When $c_{\epsilon,k}$ is small enough, $\tilde{X}$ is a symplectic hypersurface in $M$ with at most isolated singularities whose number (counting multiplicity) is uniformly bounded in each ball of radius $1/(\epsilon k)$. For suitable local complex structure, the singularities are isolated holomorphic singularities with bounded multiplicity. Furthermore, there exists $c>0$ such that $\tilde{X}$ is in $B_{\frac{c}{k\epsilon}} (\Gamma)$, and $\tilde{X}$ is smooth of bounded geometry at the $\frac{1}{k \epsilon}$-scale outside of $B_{\frac{c}{k\epsilon}} (\Gamma^{(n-1)})$. In particular, singularities of $\tilde{X}$ is in $B_{\frac{c}{k\epsilon}} (\Gamma^{(n-1)})$.
\end{theorem}
{\bf Proof:} Since $\tilde{X}$ coincides with $X$ outside of $\tilde{B}_{3R_1} (\gamma^\circ_{\tilde{I}})$, where $|d\mu_p (p)| \geq C_3$ according to proposition \ref{ah}, and $\tilde{B}_{3R_1} (\gamma^\circ_{\tilde{I}})$ is in $B_{\frac{c}{k\epsilon}} (\Gamma)$ for $c>0$ suitably large, by theorem \ref{dc}, we only need to consider the singularities of $\tilde{X}$ in $\tilde{B}_{3R_1} (\gamma^\circ_{q'_i})$ for each $i \in \tilde{I}$.\\

For $p \in \tilde{B}_{R_1, 3R_1} (\gamma^\circ_{q'_i})$, $|\tilde{\mu}_p (p)|_{C^1} \geq C_3 -Cc_{\epsilon,k}$. For $c_{\epsilon,k}$ small, $\tilde{X} \cap \tilde{B}_{R_1, 3R_1} (\gamma^\circ_{q'_i})$ is smooth of bounded geometry at the $\frac{1}{k \epsilon}$-scale.\\

On $\tilde{B}_{R_1} (\gamma^\circ_{q'_i})$, $\tilde{\mu}_{q_i} (z) = e^{m_* \cdot z} \mu^\circ_{q_i} (z)$. Hence $\tilde{X} \cap \tilde{B}_{R_1} (\gamma^\circ_{q'_i})$ has only bounded number of isolated holomorphic singularities with bounded multiplicity according to proposition \ref{ac}.
\hfill$\Box$\\
\begin{theorem}
\label{da}
When $c_{\epsilon,k}$ is small enough, $\hat{X}$ is a smooth symplectic hypersurface in $M$ of locally bounded geometry at the $\frac{1}{k \epsilon}$-scale. Furthermore, there exists $c>0$ such that $\tilde{X}$ is in $B_{\frac{c}{k\epsilon}} (\Gamma)$.
\end{theorem}
{\bf Proof:} Since $\hat{X}$ coincides with $\tilde{X}$ outside of $\tilde{B}_{2R_1} (\gamma^\circ_{\tilde{I}})$, by theorem \ref{dd}, we only need to consider $\hat{X}$ in $\tilde{B}_{2R_1} (\gamma^\circ_{q'_i})$ for each $i \in \tilde{I}$.\\

For $p \in \tilde{B}_{R_1, 2R_1} (\gamma^\circ_{q'_i})$, $|\hat{\mu}_p (p)|_{C^1} \geq C_3 -Cc_{\epsilon,k} - C'\hat{\epsilon}_i$. For $c_{\epsilon,k}$ and $\hat{\epsilon}_i$ small, $\hat{X} \cap \tilde{B}_{R_1, 2R_1} (\gamma^\circ_{q'_i})$ is smooth of bounded geometry at the $\frac{1}{k \epsilon}$-scale.\\

On $\tilde{B}_{R_1} (\gamma^\circ_{q'_i})$, $\hat{\mu}_{q_i} (z)= e^{m_* \cdot z} (\mu^\circ_{q_i} (z) + \hat{\epsilon}_i)$. We can choose $\hat{\epsilon}_i$ suitably so that there exists $C_4>0$ such that $|\hat{\mu}_p (p)|_{C^1} \geq C_4$ on $\tilde{B}_{R_1} (\gamma^\circ_{q'_i})$. Hence $\hat{X} \cap \tilde{B}_{R_1} (\gamma^\circ_{q'_i})$ is smooth of bounded geometry at the $\frac{1}{k \epsilon}$-scale.
\hfill$\Box$\\

{\bf Remark:} Although $\tilde{X}$ is not smooth, it has very mild analytic singularities and is more canonical than the smooth $\hat{X}$, therefore, it could be more useful for some applications. (For example, the artificial smoothing $\hat{X}$ hides the rather canonical singularities of $\tilde{X}$ (related to irrationality of $\{p_i\}_{i\in I}$) that usually can not be avoided by perturbation of $\{p_i\}_{i\in I}$.) Theorem \ref{da} can be viewed as a refined version of Donaldson's main theorem in \cite{D1}, here the structure of the symplectic hypersurfaces can be more explicitly known through real skeletons and basic exponential sums.\\
\begin{co}
\label{daa}
$\hat{X}$ (resp. $\tilde{X}$) represents a current on $M$.
\end{co}
{\bf Proof:} Since $\hat{X}$ is smooth, it obviously represents a current on $M$. Since the singularities of $\tilde{X}$ are analytic isolated singularities with respect to certain local complex structures of $M$, and analytic varieties represents currents, $\tilde{X}$ also represents a current on $M$.
\hfill$\Box$\\

{\bf Remark:} It is not clear at all if $X$ represents a current or not.\\

\stepcounter{subsection}
{\bf \S \thesubsection\ The symplectic Lefschetz pencil theorem:} We can also prove an alternative version of Donaldson's main theorem in \cite{D2} concerning pencil in almost the same way as the proof of theorem \ref{da}. As we observed in the remark at the beginning of section 3, we only need to require the exponential sums to be basic in the pencil case. This observation in fact enable us to give a slightly more precise proof for the pencil case in the following.\\

According to the proof of proposition \ref{af}, $I$ is separated into $N$ groups $I_j$ for $1\leq j \leq N$ so that $p_i$ for $i$ in the same group are not adjacent to each other. Let $\zeta_j = e^{\frac{2j\pi i}{N}}$ and $a_{0,i} = -a_{\infty,i} \zeta_j$ for $i \in I_j$. Let $s_t = s_0 + ts_\infty$, where $s_0 =\underset{j\in I}{\sum}a_{0,j} \sigma_{p_j}$, $s_\infty = \underset{i\in I}{\sum} a_{\infty,i} \sigma_{p_i}$; $|a_{0,i}| = |a_{\infty,i}| =1$ for $i\in I$.\\

{\bf Remark:} The same arguments will also apply when $\log |a_{0,i}| = \log |a_{\infty,i}|$ are uniformly bounded for $i\in I$. In fact, the same arguments with straightforward modification can also deal with the case when $|a_{0,i}|$ is not necessarily equals to $|a_{\infty,i}|$, as long as $a_{0,i}/a_{\infty,i}$ are well separated in $\mathbb{CP}^1$ for adjacent $p_i$'s. In such situation, $\Upsilon$ is slightly more complicated, which results in suitable modification of the arguments accordingly.\\

For $1 \leq j \leq N$, let $\Gamma_j$ be the real skeleton determined by $\{p_i\}_{i\in I\setminus I_j}$. (See the second picture of figure 2 in page \pageref{de}.) For $n \leq l \leq 2n-1$, let $\tilde{\Gamma}^{(l)}$ be the closure of {\footnotesize $\displaystyle \bigcup_{j=1}^N \bigcup_{i \in I_j}$}$(\Gamma_j^{(l)} \cap U_i^\circ)$, and $\tilde{\Gamma} = \tilde{\Gamma}^{(2n-1)}$. (See the first picture of figure 2 in page \pageref{de}.) For a tree $\Upsilon$ with a unique $N$-valent vertex $\tau^\circ$ and $N$ legs $\Upsilon_j \cong [0, -\infty)$ attach to $\tau^\circ$ for $1 \leq j \leq N$, $\{\tilde{\alpha}_\tau \}_{\tau\in \Upsilon} =$ {\footnotesize $\displaystyle \bigcup_{j=1}^N$}$\{\tilde{\alpha}_\tau \}_{\tau\in \Upsilon_j}$ is a tree in $\mathbb{R}^{N}$, where $\{\tilde{\alpha}_\tau \}_{\tau\in \Upsilon_j} = [0, -\infty)\underset{i\in I_j}{\sum} e_i$ is a ray starting from $\tilde{\alpha}_{\tau^\circ} =0$. Let $\tilde{\Gamma}_\tau$ be the real skeleton of $\underset{i\in I}{\sum}e^{\tilde{\alpha}_{\tau,i}} \sigma_{p_i}$. (See the third picture of figure 2 in page \pageref{de}.) $\{\tilde{\Gamma}_\tau\}_{\tau \in \Upsilon}$ can be understood as a \textsf{real pencil} of real hypercomplices in $M$ (parameterized by the tree $\Upsilon$). We call the real pencil $\{\tilde{\Gamma}_\tau\}_{\tau \in \Upsilon}$ the \textsf{skeleton} of the complex pencil $\{X_t\}$.\\

\pspicture(-6,-2.5)(6,1.7)
\rput(-4,0){
\psset{linecolor=lightgray}
\gray
\psline(.7,1)(1.3,1.65)
\psline(-.74,.95)(-1.3,1.6)
\psline(1.02,-.4)(1.85,-.65)
\psline(-1,-.5)(-1.8,-.7)
\psline(0,-1.2)(0,-1.7)
\qdisk(0,0){1.5pt}
\rput(.3,.2){$p_i$}
\rput(-.4,.2){$U_i$}
\qdisk(0,2){1pt}
\rput(.3,1.8){$p_{i_1}$}
\rput(-.1,1.3){$U_{i_1}$}
\qdisk(-1.8,.3){1pt}
\rput(-1.6,.6){$p_{i_2}$}
\rput(-1.3,0){$U_{i_2}$}
\qdisk(-1.3,-1.5){1pt}
\rput(-1.5,-1.3){$p_{i_3}$}
\rput(-.6,-1.2){$U_{i_3}$}
\qdisk(1.3,-1.5){1pt}
\rput(1.6,-1.3){$p_{i_4}$}
\rput(.7,-1.2){$U_{i_4}$}
\qdisk(1.8,.4){1pt}
\rput(2.1,.2){$p_{i_5}$}
\rput(1.3,.5){$U_{i_5}$}
\rput(-1.9,-.7){$\Gamma$}
\pspolygon(.7,1)(-.74,.95)(-1,-.5)(0,-1.2)(1.02,-.4)
\psset{linecolor=black}
\black
\psline(-.05,.2)(.7,1)
\psline(-.05,.2)(-.74,.95)
\psline(0,-.05)(1.02,-.4)
\psline(-.05,-.2)(-1,-.5)
\psline(-.05,-.2)(0,-1.2)
\psline(-.05,.2)(0,-.05)
\psline(-.05,-.2)(0,-.05)
\psline(.7,1)(1.6,1.05)
\psline(.7,1)(.48,1.9)
\psline(-.74,.95)(-1.64,.9)
\psline(-.74,.95)(-.60,1.75)
\psline(1.02,-.4)(1.24,-1.3)
\psline(1.02,-.4)(1.72,.1)
\psline(-1,-.5)(-1.15,-1.3)
\psline(-1,-.5)(-1.7,0)
\psline(0,-1.2)(.7,-1.7)
\psline(0,-1.2)(-.7,-1.7)
\rput(1.8,1.1){$\tilde{\Gamma}$}}

\pspolygon[linecolor=lightgray](.7,1)(-.74,.95)(-1,-.5)(0,-1.2)(1.02,-.4)
\psline(-.05,.2)(1.3,1.65)
\psline(-.05,.2)(-1.3,1.6)
\psline(0,-.05)(1.85,-.65)
\psline(-.05,-.2)(-1.8,-.7)
\psline(-.05,-.2)(0,-1.7)
\psline(-.05,.2)(0,-.05)
\psline(-.05,-.2)(0,-.05)
\rput(1.3,-.75){$\Gamma_j$}

\rput(4,0){
\psline[linecolor=lightgray](-.05,.2)(1.3,1.65)
\psline[linecolor=lightgray](-.05,.2)(-1.3,1.6)
\psline[linecolor=lightgray](0,-.05)(1.85,-.65)
\psline[linecolor=lightgray](-.05,-.2)(-1.8,-.7)
\psline[linecolor=lightgray](-.05,-.2)(0,-1.7)
\psline[linecolor=lightgray](-.05,.2)(0,-.05)
\psline[linecolor=lightgray](-.05,-.2)(0,-.05)
\pspolygon[linecolor=lightgray](.7,1)(-.74,.95)(-1,-.5)(0,-1.2)(1.02,-.4)
\psline(.4,.7)(1.3,1.65)
\psline(-.44,.65)(-1.3,1.6)
\psline(.6,-.25)(1.85,-.65)
\psline(-.6,-.37)(-1.8,-.7)
\psline(-.03,-.78)(0,-1.7)
\rput(1.3,-.75){$\tilde{\Gamma}_\tau$}
\pspolygon(.4,.7)(-.44,.65)(-.6,-.37)(-.03,-.78)(.6,-.25)}
\stepcounter{figure}
\label{de}
\uput{2}[d](0,0){Figure \thefigure: The real skeletons $\tilde{\Gamma}$, $\Gamma_j$ ($i\in I_j$) and $\tilde{\Gamma}_\tau$ ($\tau \in \Upsilon_j$)}
\endpspicture

Let $D_j = \{t \in \mathbb{C}: |\log (t/\zeta_j)| \leq \frac{\pi}{N}\}$ and $D'$ be the compliment of $\{D_j\}_{j=1}^N$ in $\mathbb{CP}^1$. Define $\tau_t = \log |\frac{N}{\pi} \log (t/\zeta_j)| \in \Upsilon_j$ for $t\in D_j$ and $\tau_t =\tau^\circ$ for $t\in D'$. $t\rightarrow \tau_t$ defines a continuous and surjective map $\mathbb{CP}^1 \setminus \{\zeta_j\}_{j=1}^N \rightarrow \Upsilon$. We will use $\tau_{\zeta_j}$ to denote the infinity of $\Upsilon_j$. $\tilde{\Gamma}_{\tau_{\zeta_j}}$ can be identified with $\Gamma_j$.\\
\begin{prop}
\label{ed}
{\footnotesize $\displaystyle \bigcap_{\tau\in \Upsilon}$}$\tilde{\Gamma}_\tau^{(l)} = \Gamma^{(l-1)} = \tilde{\Gamma}^{(l)} \cap \Gamma$ for $n\leq l\leq 2n-1$. In particular, the base locus of the real pencil is {\footnotesize $\displaystyle \bigcap_{\tau\in \Upsilon}$}$ \tilde{\Gamma}_\tau = \Gamma^{(2n-2)} = \tilde{\Gamma} \cap \Gamma$. {\footnotesize $\displaystyle \bigcup_{\tau\in \Upsilon}$}$\tilde{\Gamma}_\tau^{(l-1)} \subset \tilde{\Gamma}^{(l)}$ for $n\leq l\leq 2n-1$, and equality hold when $\{p_i\}_{i\in I}$ is generic (as in proposition \ref{af}).
\end{prop}
{\bf Proof:} For $n\leq l \leq 2n-1$, when $\tau \in \Upsilon_j$ approaches infinity, all the top dimension strata of $\Gamma^{(l)}$ will expand in $\tilde{\Gamma}_\tau^{(l)}$ except those in $U_i$ for $i\in I_j$, which will move and shrink in $\tilde{\Gamma}_\tau^{(l)}$. (See figure 2.) Since every strata of $\Gamma^{(l-1)}$ belongs to a top dimension strata of $\Gamma^{(l)}$ that is not in $U_i$ for all $i\in I_j$, we have $\Gamma^{(l-1)} \subset \tilde{\Gamma}_\tau^{(l)}$. On the other hand, every top dimension strata of $\Gamma^{(l)}$ is in some $U_i$ for $i$ in some $I_j$, and will not be in $\tilde{\Gamma}_\tau^{(l)}$ when $\tau \in \Upsilon_j$ approaches infinity. Hence {\footnotesize $\displaystyle \bigcap_{\tau\in \Upsilon}$}$ \tilde{\Gamma}_\tau^{(l)} = \Gamma^{(l-1)}$.\\

Since $\Gamma_j^{(l)} \cap \partial U_i = (\partial U_i)^{(l-1)}$ for $i\in I_j$, by definition of $\tilde{\Gamma}$, $\tilde{\Gamma}^{(l)} \cap \Gamma = \Gamma^{(l-1)}$ when $n\leq l \leq 2n-1$.\\

When $\tau \in \Upsilon_j$, $(\tilde{\Gamma}_\tau \setminus \Gamma_j) \cap U_i$ for $i\in I_j$ only contains top dimensional open strata of $\tilde{\Gamma}_\tau$. Hence $\tilde{\Gamma}_\tau^{(l-1)} \cap U_i \subset \Gamma_j^{(l)} \cap U_i = \tilde{\Gamma}^{(l)} \cap U_i$ for $n\leq l\leq 2n-1$. Outside of $U_i$ for $i\in I_j$, $\tilde{\Gamma}_\tau^{(l-1)} = \Gamma^{(l-1)}$. Consequently, $\tilde{\Gamma}_\tau^{(l-1)} \subset \tilde{\Gamma}^{(l)}$ for $\tau \in \Upsilon_j$.\\

For the equality to fail, there should be a $l$-strata of $\tilde{\Gamma}$ in the interior of $U_i$ for $i$ in some $I_j$, which is in the limit of $\tilde{\Gamma}_\tau \setminus \Gamma_j$ when $\tau \in \Upsilon_j$ approaches infinity. This does not happen for generic $\{p_i\}_{i\in I}$, in which case, the diameter of $(\tilde{\Gamma}_\tau \setminus \Gamma_j) \cap U_i$ approaches $0$ when $\tau \in \Upsilon_j$ approaches infinity.
\hfill$\Box$\\
\begin{prop}
\label{aj}
For $p \in X_t$, let $s_t = \mu_{t,p} \sigma_p$. We may write $\mu_{t,p} = \mu^\circ_{t,p} + \mu'_{t,p}$ so that

\[
\mu^\circ_{t,p} = \mbox{{\small $\displaystyle\sum_{i\in I_p}$ }} a_{t,i} e^{-\epsilon^2 k|m_i|^2} e^{m_i \cdot z}
\]

is a basic exponential sum. Furthermore, $|\mu'_{t,p}|_{C^1_p} = O(c_{\epsilon,k})$ when $|z|$ is bounded, where $c_{\epsilon,k} = \max(\epsilon, 1/(\epsilon k), e^{-c\epsilon^2k})$ for certain fixed $c>0$.
\end{prop}
{\bf Proof:} Since $|a_{0,i}| = |a_{\infty,i}|=1$, we have $a_{t,i} = O(1+|t|)$. This proposition is then a direct consequence of proposition \ref{ag} and its proof if one can show that $\max_{i\in I_p} (|a_{t,i}|) \sim 1+|t|$ uniformly for $t\in \mathbb{CP}^1$.\\

Proposition \ref{af} implies that $\mu^\circ_{t,p}$ is a pencil of strongly basic exponential sums. According to proposition \ref{bb} and the special structure of the tree $\Upsilon$ therein, for all (except possibly one) $i\in I_p$, we have $|a_{t,i}| \sim 1+|t|$ uniformly for $t\in \mathbb{CP}^1$. Our definition of $I_p$ clearly implies that $|I_p|\geq 2$. Consequently, $\max_{i\in I_p} (|a_{t,i}|) \sim 1+|t|$ uniformly for $t\in \mathbb{CP}^1$.
\hfill$\Box$\\
\begin{prop}
For $c_{\epsilon,k}$ small, there exists $c>0$ such that $X_t \subset B_{\frac{c}{\epsilon k}} (\tilde{\Gamma}_{\tau_t})$ for any $t\in \mathbb{CP}^1$.
\end{prop}
{\bf Proof:} The proposition is a consequence of propositions \ref{aj}, \ref{bb} and the definition of $\tau_t$.
\hfill$\Box$\\
\begin{prop}
\label{ah}
For $c_{\epsilon,k}$ small, there exist $R_1,C_3>0$, a set of points $\{q'_i\}_{i\in \tilde{I}} \subset M$ and $\{\gamma^\circ_{q'_i}\}_{i\in \tilde{I}}$, where $\gamma^\circ_{q'_i}$ is a set of singular points of the pencil $\{\mu^\circ_{t,p}\}$, so that $\{\tilde{B}_{3R_1} (\gamma^\circ_{q'_i})\}_{i\in \tilde{I}}$ are disjoint and $|\mu_{t,p}|_{C^1_p} (p) \geq C_3$ for $p \in M \setminus \tilde{B}_{R_1} (\gamma^\circ_{\tilde{I}})$, where $\gamma^\circ_{\tilde{I}} = \bigcup_{i \in \tilde{I}} \gamma^\circ_{q'_i}$.
\end{prop}
{\bf Proof:} Let $R_1$ satisfy $14R_1 N_{B_1} <1$. Proposition \ref{bd} implies that there exists $C_1>0$ such that $\displaystyle \min_{t \in \mathbb{CP}^1} \left( |\mu^\circ_{t,p}|_{C^1_p} \right) \geq C_1$ away from $R_1$-balls centered at singular points of the pencil $\{\mu^\circ_{t,p}\}$. Proposition \ref{ag} implies that for $c_{\epsilon,k}$ small, there exists a constant $C_2>0$ such that $|\mu_{t,p}|_{C^1_p} (z) \leq C_2$ implies that $|\mu^\circ_{t,p}|_{C^1_p} (z) \leq C_1$.\\

Let $q'_0 \in M$ be a point satisfying $\displaystyle \min_{t \in \mathbb{CP}^1} \left(|\mu_{t,q'_0}|_{C^1_{q'_0}} (q'_0)\right) < C_3 \leq C_2$, then $\displaystyle \min_{t \in \mathbb{CP}^1} \left(|\mu^\circ_{t,q'_0}|_{C^1_{q'_0}} (q'_0)\right) < C_1$. $q'_0$ is in the $R_1$-ball of a singular point $q_0$ of the pencil $\{\mu^\circ_{t,q'_0}\}$. Let $\gamma_{q'_0}$ be the maximal connected graph in $\tilde{B}_1(q'_0)$ with the vertex set $\gamma^\circ_{q'_0} = \{q_{0,j}\}_{j\in A_0}$ (including $q_{0,0} = q_0$) containing singular points of the pencil $\{\mu^\circ_{t,p}\}$ such that any leg of $\gamma_{q'_0}$ has length less than $7R_1$. The condition $14R_1 N_{B_1} <1$ implies that $\gamma_{q'_0} \subset \tilde{B}_{1/2} (q'_0)$.\\

Take $q'_1 \in M \setminus \tilde{B}_{R_1} (\gamma^\circ_{q'_0})$ satisfying $\displaystyle \min_{t \in \mathbb{CP}^1} \left(|\mu_{t,q'_1} |_{C^1_{q'_1}} (q'_1)\right) < C_3 \leq C_2$, by similar procedure we may construct $\gamma_{q_1}$ with vertices $\{q_{1,j}\}_{j\in A_1}$ (including $q_{1,0} = q_1$) that satisfies the additional condition: $q_{1,j} \not\in \tilde{B}_{R_1} (\gamma^\circ_{q'_0})$ for $j\in A_1$. For $j_1 \in A_1$, $\displaystyle \min_{t \in \mathbb{CP}^1} \left(|\mu^\circ_{t,q'_1} |_{C^1_{q'_1}} (q_{1,j_1})\right) =0$ implies that $\displaystyle \min_{t \in \mathbb{CP}^1} \left(|\mu_{t,q'_1} |_{C^1_{q'_1}} (q_{1,j_1})\right) < C_3$. By proposition \ref{ai}, when $C_3$ is small, $\displaystyle \min_{t \in \mathbb{CP}^1} \left(|\mu_{t,q'_0} |_{C^1_{q'_0}} (q_{1,j_1})\right) < C_2$,  $\displaystyle \min_{t \in \mathbb{CP}^1} \left(|\mu^\circ_{t,q'_0}|_{C^1_{q'_0}} (q_{1,j_1})\right) < C_1$. By proposition \ref{bd}, $q_{1,j_1}$ is in the $R_1$-ball of a singular point $q_{0,j_0}$ of the pencil $\{\mu^\circ_{t,q'_0}\}$, which together with $q_{1,j_1} \not\in \tilde{B}_{R_1} (\gamma^\circ_{q'_0})$ implies that $j_0 \not\in A_0$. Hence  $q_{0,j_0} \not\in \tilde{B}_{7R_1} (\gamma^\circ_{q'_0})$,  $q_{1,j_1} \not\in \tilde{B}_{6R_1} (\gamma^\circ_{q'_0})$. Consequently, $\tilde{B}_{3R_1} (\gamma^\circ_{q'_0})$ is disjoint from $\tilde{B}_{3R_1} (\gamma^\circ_{q'_1})$.\\

It remains to verify that $q_{1,0} = q_1 \not\in \tilde{B}_{R_1} (\gamma^\circ_{q'_0})$. Since $\displaystyle \min_{t \in \mathbb{CP}^1} \left(|\mu_{t,q'_1} |_{C^1_{q'_1}} (q'_1)\right) \leq C_3$, by proposition \ref{ai}, when $C_3$ is small, $\displaystyle \min_{t \in \mathbb{CP}^1} \left(|\mu_{t,q'_0} |_{C^1_{q'_0}} (q'_1)\right) \leq C_2$, and $\displaystyle \min_{t \in \mathbb{CP}^1} \left(|\mu^\circ_{t,q'_0} |_{C^1_{q'_0}} (q'_1)\right) \leq C_1$. By proposition \ref{ac}, $q'_1$ is in the $R_1$-ball of a critical point $q_{0,j_0}$ of $\mu^\circ_{q'_0}$, which together with $q'_1 \not\in \tilde{B}_{R_1} (\gamma^\circ_{q'_0})$ imply that $j_0 \not\in A_0$. Hence $q_{0,j_0} \not\in \tilde{B}_{7R_1} (\gamma^\circ_{q'_0})$, which together with $q_{1,0} \in \tilde{B}_{2R_1} (q_{0,j_0})$ imply that $q_{1,0} \not\in \tilde{B}_{2R_1} (\gamma^\circ_{q'_0})$.\\

Through induction, we can construct $\{\gamma_{q_i}\}_{i\in \tilde{I}}$ so that $\displaystyle \min_{t \in \mathbb{CP}^1} \left(|\mu_{t,p}|_{C^1_p} (p)\right) \geq C_3$ for $p \in X \setminus \tilde{B}_{R_1} (\gamma^\circ_{\tilde{I}})$ and $\{\tilde{B}_{3R_1} (\gamma^\circ_{q'_i})\}_{i\in \tilde{I}}$ are disjoint.
\hfill$\Box$\\

For each $i\in \tilde{I}$, we can construct a cutoff function $\tilde{\rho}_i$ (resp. $\hat{\rho}_i$) such that $\tilde{\rho}_i =1$ (resp. $\hat{\rho}_i =1$) on $\tilde{B}_{2R_1} (\gamma^\circ_{q'_i})$ (resp. $\tilde{B}_{R_1} (\gamma^\circ_{q'_i})$) and $\tilde{\rho}_i =0$ (resp. $\hat{\rho}_i =0$) away from $\tilde{B}_{3R_1} (\gamma^\circ_{q'_i})$ (resp. $\tilde{B}_{2R_1} (\gamma^\circ_{q'_i})$). Let $\hat{X}_t = \{\hat{s}_t=0\}$ and $\tilde{X}_t = \{\tilde{s}=0\}$, where

\[
\tilde{s}_t = s_t - \mbox{{\small $\displaystyle\sum_{i\in \tilde{I}}$ }} \tilde{\rho}_i \mu'_{t,q'_i} \sigma_{q'_i}, \ \hat{s}_t = \tilde{s}_t + \mbox{{\small $\displaystyle\sum_{i\in \tilde{I}}$ }} \hat{\epsilon}_{i,t} \hat{\rho}_i \sigma_{q'_i},
\]

$\hat{\epsilon}_{i,t} = \hat{\epsilon}_{i,0} + t\hat{\epsilon}_{i,\infty}$, $\hat{\epsilon}_{i,0}$ and $\hat{\epsilon}_{i,\infty}$ are suitable small constants such that $0 \in \mathbb{C}^2$ is not a critical value of $(\mu^\circ_{0,q'_i} + \hat{\epsilon}_{i,0}, \mu^\circ_{\infty,q'_i} + \hat{\epsilon}_{i,\infty})$ on $\tilde{B}_{R_1} (\gamma_{q'_i})$.\\

\begin{theorem}
\label{db}
When $c_{\epsilon,k}$ is small enough, $\tilde{X}_t$ (resp. $\hat{X}_t$) is a pencil of symplectic hypersurface in $M$ with at most isolated singularities whose number (counting multiplicity) is uniformly bounded in each ball of radius $1/(\epsilon k)$. For suitable local complex structure, the singularities are isolated holomorphic singularities with bounded multiplicity. $\hat{X}_t$ is a \textsf{generalized Lefschetz pencil} of symplectic hypersurfaces in $M$, in the sense that in addition, the base locus $\hat{Y}$ is smooth, which implies the bounded finiteness of the number of singular fibres. Furthermore, There exists $c>0$ such that the singular set of the pencil $\{X_t\}$ (resp. $\{\hat{X}_t\}$ and $\{\tilde{X}_t\}$) is in $B_{\frac{c}{k\epsilon}} (\tilde{\Gamma}^{(n)})$, and the singular set of the base locus $Y$ (resp. $\tilde{Y}$) for $\{X_t\}$ (resp. $\{\tilde{X}_t\}$) is in $B_{\frac{c}{k\epsilon}} (\Gamma^{(n-1)})$.
\end{theorem}
{\bf Proof:} Since $\hat{X}_t$ (resp. $\tilde{X}_t$) coincide with $X_t$ outside of $\tilde{B}_{3R_1} (\gamma^\circ_{\tilde{I}})$, where $|\mu_{t,p}|_{C^1_p} (p) \geq C_3$, and $\tilde{B}_{3R_1} (\gamma^\circ_{\tilde{I}})$ is in $B_{\frac{c}{k\epsilon}} (\Gamma)$ for $c>0$ suitably large, we only need to consider $\tilde{B}_{3R_1} (\gamma^\circ_{q'_i})$ for each $i \in \tilde{I}$.\\

For $p \in \tilde{B}_{R_1, 3R_1} (\gamma^\circ_{q'_i})$, $|\tilde{\mu}_{t,p}|_{C_p^1} (p), |\hat{\mu}_{t,p}|_{C_p^1} (p) \geq C_3 -Cc_{\epsilon,k} - C'\hat{\epsilon}_i$. For $c_{\epsilon,k}$ and $\hat{\epsilon}_i$ small, $\hat{X}_t \cap \tilde{B}_{2R_1, 3R_1} (\gamma^\circ_{q'_i})$ and $\tilde{X}_t \cap \tilde{B}_{2R_1, 3R_1} (\gamma^\circ_{q'_i})$ are smooth.\\

On $\tilde{B}_{R_1} (\gamma^\circ_{q'_i})$, $\tilde{\mu}_{t,q_i} (z) = \mu^\circ_{t,q_i} (z)$, $\hat{\mu}_{t,q_i} (z)= \mu^\circ_{t,q_i} (z) + \hat{\epsilon}_{i,t}$. Hence $\hat{X} \cap \tilde{B}_{R_1} (\gamma^\circ_{q'_i})$ (resp. $\tilde{X} \cap \tilde{B}_{R_1} (\gamma^\circ_{q'_i})$) has only isolated holomorphic singularities with bounded multiplicity and $\hat{Y} \cap \tilde{B}_{R_1} (\gamma^\circ_{q'_i})$ is smooth.
\hfill$\Box$\\

{\bf Remark:} It is straightforward to get a Lefschetz pencil by perturbing the finite many isolated singular points (with bounded multiplicity) of the pencil $\{\hat{X}_t\}$ to non-degenerate singular points. Similar to the hypersurface case, the pencils $\{\hat{X}_t\}$ and $\{\tilde{X}_t\}$ are more canonical and structured, such further perturbation is not essential because isolated singularity with bounded multiplicity $l$ is well understood and can be easily perturb to $l$ non-degenerate singularities locally.\\

\se{The limit of currents}
According to corollary \ref{daa}, $\tilde{X}$ and $\hat{X}$ both represent $(2n-2)$-currents denoted as $[\tilde{X}]$ and $[\hat{X}]$. In this section, we examine the limit of currents $\frac{1}{k}[\tilde{X}]$ and $\frac{1}{k}[\hat{X}]$ as $k\rightarrow +\infty$. In Donaldson's case (\cite{D1}), such limit of currents is $\omega$ as a $(2n-2)$-current, which indicates that the symplectic hypersurfaces in \cite{D1} distribute quite evenly through out $M$ for $k$ large. In our case, the symplectic hypersurfaces are concentrated near the real skeletons and the limiting current is determined by the limit of the real skeletons. Assume $D\tilde{s} = \tilde{A} \tilde{s}$ (resp. $D\hat{s} = \hat{A} \hat{s}$), then $-ik\omega = d\tilde{A} + [\tilde{X}]$ (resp. $-ik\omega = d\hat{A} + [\hat{X}]$).\\

For $i\in I$, recall that $-i\omega = \frac{1}{2} dJ_{p_i}d \log h_i$ on $U_i$, where $h_i = e^{-|p-p_i|_{p_i}^2/2}$. Let $A$ denote the connection 1-form of $L^k$ with discontinuous coefficients such that $A = A_i = \frac{k}{2} J_{p_i}d \log h_i$ on $U_i$. Then we have the equation  $-i\omega = \frac{1}{k} d A + [\beta_{\Gamma}]$ of currents. $\beta_{\Gamma}$ can be viewed as a smooth 1-form defined on the smooth part $\Gamma^\circ$ of $\Gamma$. Let $\Gamma_{ij} = \{p\in M : h_i(p) = h_j(p)\}$. When $U_j$ is adjacent to $U_i$, $\Gamma_{ij}$ contains part of $\Gamma$ (that is $U_i \cap U_j$), where $\beta_{\Gamma} = \frac{1}{2} Jd(\log h_i - \log h_j)$. For any smooth $(2n-2)$-form $\psi$ on $M$,

\[
\langle [\tilde{X}], \psi \rangle = \int_{\tilde{X}} \psi,\ \langle [\hat{X}], \psi \rangle = \int_{\hat{X}} \psi,\ \langle [\beta_{\Gamma}], \psi \rangle = \int_{\Gamma^\circ} \beta_{\Gamma} \wedge \psi.
\]

Let $\displaystyle M \setminus B_{\frac{c}{\epsilon k}}(\Gamma) = \bigcup_{i\in I} U_{i,\epsilon}$, $B^\circ_{\frac{c}{\epsilon k}}(\Gamma) = B_{\frac{c}{\epsilon k}}(\Gamma) \setminus \tilde{B}_{R_1} (\gamma^\circ_{\tilde{I}})$, where $U_{i,\epsilon} \subset U_i$ is a neighborhood of $p_i$. $B^\circ_{\frac{c}{\epsilon k}}(\Gamma)$, $\{U_{i,\epsilon}\}_{i\in I}$ and $\{\tilde{B}_{3R_1} (\gamma^\circ_{q'_i})\}_{i\in \tilde{I}}$ form a cover of $M$.\\
\begin{prop}
\label{ea}
(1) For a complex valued $C^\infty$-function $f$ that is transverse to $0$, $\displaystyle \int_{B(p)} \frac{1}{|f|}$ is finite (with the bound depending on the positive lower bound of $|df|$ near the zero set of $f$).\\

(2) For any holomorphic function $f$, $\displaystyle \int_{B(p)} \log |f|$ is finite. For a basic exponential sum $\mu$, $\displaystyle \int_{B(p)} \log |\mu|_p$ is finite (with the bound depending on the geometry of the set of exponents of $\mu$).\\
\end{prop}
{\bf Proof:} (1) is a simple exercise that come down to the fact: $\int_{B_1} \frac{1}{|z|} \leq C$ for $B_1 \subset \mathbb{C}$. The first part of (2) is a somewhat non-trivial well known fact in several complex variables, whose proof uses Weierstrass preparation theorem.\\

For the second part of (2), we will use the by now familiar limiting method. Assume there is a sequence $(\mu_i, q_i)$ such that

\[
\lim_{i \rightarrow +\infty} \int_{B(q_i)} \log |\mu_i|_{q_i} = \infty.
\]

Then let $\tilde{\mu}_i (z) = e^{-b_i(q_i)} \mu_i (q_i +z)$. By possibly taking subsequence, we may assume the limit $\displaystyle \tilde{\mu}_\infty = \lim_{i \rightarrow +\infty} \tilde{\mu}_i$ exists and is not zero. Since $\tilde{\mu}_i$ are holomorphic functions, the convergence is uniform on $B(0)$. According to Weierstrass preparation theorem, there exists $R_1,R_2>0$ and suitable coordinate $z = (w,z')$ such that $\tilde{\mu}_\infty = h(z) P(z)$, where $h(z) \not=0$ for $z\in B_{2R_1}^w \times B_{R_2}^{z'}$ and

\[
P(z) = w^l + a_1(z') w^{l-1} + \cdots + a_i(z')
\]

is a polynomial on $w$, whose roots $\{w_j(z')\}_{j=1}^l$ are all in $B_{R_1}^w$, when $z' \in B_{R_2}^{z'}$. Since $\tilde{\mu}_i$ converges to $\tilde{\mu}_\infty$ uniformly, when $i$ is large enough, the corresponding Weierstrass polynomial has the same degree $l$ with the roots $\{w_{i,j}(z')\}_{j=1}^l$ that converge to $\{w_j(z')\}_{j=1}^l$ uniformly. Let $\rho (w)$ be a smooth function supported in $B_{2R_1}^w$ such that $\Delta_w \rho = 1$ on $B_{R_1}^w$.

\[
\int_{B_{2R_1}^w \times B_{R_2}^{z'}} \log |\tilde{\mu}_\infty| \Delta_w \rho = \sum_{j=1}^l \int_{B_{R_2}^{z'}} \rho (w_j(z')).
\]
\[
\int_{B_{R_1}^w \times B_{R_2}^{z'}} \log |\tilde{\mu}_\infty| = \sum_{j=1}^l \int_{B_{R_2}^{z'}} \rho (w_j(z')) - \int_{B_{R_1,2R_1}^w \times B_{R_2}^{z'}} \log |\tilde{\mu}_\infty| \Delta_w \rho.
\]

Since $\log |\tilde{\mu}_i|$ (resp. $\{w_{i,j}(z')\}_{j=1}^l$) converges uniformly to $\log |\tilde{\mu}_\infty|$ (resp. $\{w_j(z')\}_{j=1}^l$) on $B_{R_1,2R_1}^w \times B_{R_2}^{z'}$ (resp. $B_{R_2}^{z'}$), using this formula, one can show that

\[
\int_{B_{R_1}^w \times B_{R_2}^{z'}} \log |\tilde{\mu}_\infty| = \lim_{i \rightarrow +\infty} \int_{B_{R_1}^w \times B_{R_2}^{z'}} \log |\tilde{\mu}_i|.
\]

In particular, $\displaystyle \int_{B_{R_1}^w \times B_{R_2}^{z'}} \log |\tilde{\mu}_i|$ is uniformly bounded. $B(0)$ can be covered by finitely many such $B_{R_1}^w \times B_{R_2}^{z'}$. Hence $\displaystyle \int_{B(q_i)} \log |\mu_i|_{q_i} = \int_{B(0)} \log |\tilde{\mu}_i|$ is uniformly bounded, which is a contradiction.
\hfill$\Box$\\

{\bf Remark:} (1) was used by Donaldson in \cite{D1}. We need (2) in the case of $\tilde{X}$ with possibly isolated analytic singularities.\\
\begin{lm}
\label{eb}
There exists $C>0$ such that for $p \in U_i$ and $j\not= i$,
\begin{equation}
\label{eb1}
|p-p_j|_{p_j}^2 - |p-p_i|_{p_i}^2 \geq C|p_j-p_i|_{p_i}{\rm Dist}_g(p,\Gamma_{ij}).
\end{equation}
\end{lm}
{\bf Proof:} When ${\rm Dist}_g (p_i,p_j) \geq 4\epsilon$, for $p\in U_i$, $|p-p_j|_{p_j} \geq \frac{1}{2}{\rm Dist}_g (p_i,p_j) + {\rm Dist}_g(p,\Gamma_{ij})$, $|p-p_i|_{p_i} \leq {\rm Dist}_g(p,\Gamma_{ij})$. The lemma is then obvious.\\

When $|p-p_j|_{p_j} + |p-p_i|_{p_i} = O(\epsilon)$. Let $\vec{n}$ be the unit vector along the direction $\overrightarrow{p_jp_i}$ under the coordinate $z_{p_i}$. Normal vectors of $\Gamma_{ij}$ are $\epsilon$-perturbations of $\vec{n}$.

\[
\vec{n} (|p-p_i|_{p_i}^2) = (z_{p_i} (p) - z_{p_i} (p_i), \vec{n}),\ \vec{n} (|p-p_j|_{p_j}^2) = (z_{p_i} (p) - z_{p_i} (p_j), \vec{n}) + O(\epsilon^2).
\]
\[
\vec{n} (|p-p_j|_{p_j}^2 - |p-p_i|_{p_i}^2) = (z_{p_i} (p_i) - z_{p_i} (p_j), \vec{n}) + O(\epsilon^2) = |p_j-p_i|_{p_i} + O(\epsilon^2).
\]

Together with $\vec{n} ({\rm Dist}_g(p,\Gamma_{ij})) = 1 + O(\epsilon)$, we have

\[
\vec{n} (|p-p_j|_{p_j}^2 - |p-p_i|_{p_i}^2) \geq (1 - O(\epsilon)) \vec{n} (|p_j-p_i|_{p_i}{\rm Dist}_g(p,\Gamma_{ij})).
\]

Since the estimate (\ref{eb1}) is true for $p\in \Gamma_{ij}$, we have $|p-p_j|_{p_j}^2 - |p-p_i|_{p_i}^2 \geq (1 - O(\epsilon)) |p_j-p_i|_{p_i}{\rm Dist}_g(p,\Gamma_{ij})$ for $p\in M$ in the $p_i$ side of $\Gamma_{ij}$ (including $U_i$).
\hfill$\Box$\\
\begin{theorem}
\label{eg}
For fixed $\epsilon>0$ that is small enough, $\displaystyle \lim_{k \rightarrow +\infty} \frac{1}{k} [\tilde{X}] = \lim_{k \rightarrow +\infty} \frac{1}{k} [\hat{X}] = [\beta_{\Gamma}]$ as currents. More precisely, for any smooth $(2n-2)$-form $\psi$ on $M$,
\begin{equation}
\label{eg5}
\langle [\tilde{X}], \psi \rangle - k\langle [\beta_{\Gamma}], \psi \rangle \leq \frac{C_1}{\epsilon} |\psi|_{C^1(M)} + \frac{C_2}{\epsilon^2 k} |\psi|_{C^2(M)},
\end{equation}
\begin{equation}
\label{eg6}
\langle [\hat{X}], \psi \rangle - k\langle [\beta_{\Gamma}], \psi \rangle \leq \frac{C}{\epsilon} |\psi|_{C^1(M)}.
\end{equation}
\end{theorem}
{\bf Proof:}
\[
\langle [\tilde{X}], \psi \rangle - k\langle [\beta_{\Gamma}], \psi \rangle = \int_M (A - \tilde{A}) \wedge d\psi.
\]

On $U_i \setminus \tilde{B}_{3R_1} (\gamma^\circ_{\tilde{I}})$, $\tilde{s} = s$ and

\[
\tilde{A} - A = \frac{Ds}{s} - A_i = \sum_{j\not= i}  \frac{a_j (D - A_i)\sigma_{p_j}}{s}.
\]

Recall that $A_i = \frac{k}{4}(z_{p_i} d\bar{z}_{p_i} - \bar{z}_{p_i} dz_{p_i})$ and $D \sigma_{p_j} = (A_j\rho_{p_j} + d\rho_{p_j}) e^{-k|z_{p_j}|^2/4}$.

\[
|A_i(p)|_g = O(k|p-p_i|_{p_i}),\ |D \sigma_{p_j} (p)|_g \leq Ck|p-p_j|_{p_j}e^{-k|p-p_j|_{p_j}^2/4}.
\]

\[
\frac{|a_j| (|D\sigma_{p_j}| + |A_i| |\sigma_{p_j}|)}{|a_i \sigma_{p_i}|} \leq C k|p_i-p_j|_{p_i} e^{-k(|p-p_j|_{p_j}^2 - |p-p_i|_{p_i}^2)/4}.
\]

By lemma \ref{eb},
\[
|\tilde{A} - A| \leq C\frac{|a_i \sigma_{p_i}|}{|s|} \sum_{j\not= i} k|p_i-p_j|_{p_i} e^{-k|p_i-p_j|_{p_i} {\rm Dist}_g(p,\Gamma_{ij})/2}.
\]

It is straightforward to derive that

\[
\int_{U_i} e^{-k|p_i-p_j|_{p_i}{\rm Dist}_g(p,\Gamma_{ij})/2} \leq \frac{C|U_i|}{k\epsilon |p_i-p_j|_{p_i}} e^{-k|p_i-p_j|_{p_i}{\rm Dist}_g(U_i,\Gamma_{ij})/2}.
\]

On $U_{i,\epsilon} \subset U_i \setminus \tilde{B}_{3R_1} (\gamma^\circ_{\tilde{I}})$, $C|s| \geq |a_i \sigma_{p_i}|$ and

\begin{equation}
\label{eg1}
\int_{U_{i,\epsilon}} |\tilde{A} - A| \leq \frac{C|U_i|}{\epsilon} \sum_{j\not= i} e^{-k|p_i-p_j|_{p_i}{\rm Dist}_g(U_i,\Gamma_{ij})/2} \leq \frac{C|U_i|}{\epsilon}.
\end{equation}

For $p \in B_{\frac{c}{\epsilon k}}(\Gamma) \cap U_i$, $|A_i(p)|_g = O(k|p-p_i|_{p_i}) = O(k\epsilon)$. Hence

\begin{equation}
\label{eg2}
\int_{B_{\frac{c}{\epsilon k}}(\Gamma)} |A| \leq Ck\epsilon |B_{\frac{c}{\epsilon k}}(\Gamma)| \leq \frac{C|M|}{\epsilon}.
\end{equation}

Here we are using $|B_{\frac{c}{\epsilon k}}(\Gamma)| \leq \frac{C|M|}{k\epsilon^2}$, which is a consequence of $|\Gamma| \leq \frac{C}{\epsilon}|M|$.\\

For $p \in B^\circ_{\frac{c}{\epsilon k}}(\Gamma) \cap U_i  \subset U_i \setminus \tilde{B}_{R_1} (\gamma^\circ_{\tilde{I}})$, we can similarly derive that

\[
|\tilde{A}(p)| \leq \frac{C|a_i \sigma_{p_i}(p)|}{|s(p)|} \sum_{j\in I} k|p_i-p_j|_{p_i} e^{-k|p_i-p_j|_{p_i}{\rm Dist}(p,\Gamma_{ij})/2} \leq \frac{Ck\epsilon}{|s(p)|_p}.
\]

(Strictly speaking, in $\tilde{B}_{R_1,3R_1} (\gamma^\circ_{q'_i})$, one also need to consider $\tilde{\rho}_i \mu'_{q'_i} \sigma_{q'_i}$. Since $|d\tilde{\rho}_i| = O(k\epsilon)$, this term of $\tilde{s}$ does not affect our estimate.) By (1) of proposition \ref{ea}, we have

\begin{equation}
\label{eg3}
\int_{B^\circ_{\frac{c}{\epsilon k}}(\Gamma)} |\tilde{A}| \leq Ck\epsilon |B_{\frac{c}{\epsilon k}} (\Gamma)| \leq \frac{C|M|}{\epsilon}.
\end{equation}

For $i\in \tilde{I}$, $\tilde{A} = Jd\log |\tilde{s}|$ in $\tilde{B}_{2R_1} (\gamma^\circ_{q'_i})$.

\[
\int_{\tilde{B}_{2R_1} (\gamma^\circ_{q'_i})} \tilde{A} \wedge \hat{\rho}_i d\psi = \int_{\tilde{B}_{2R_1} (\gamma^\circ_{q'_i})} \log |\tilde{s}| \wedge d (\hat{\rho}_i *^{-1}J *d\psi).
\]

From this equation, apply (2) of proposition \ref{ea}, one can derive

\[
\int_{\tilde{B}_{2R_1} (\gamma^\circ_{q'_i})} \tilde{A} \wedge \hat{\rho}_i d\psi \leq |\tilde{B}_{2R_1} (\gamma^\circ_{q'_i})| (C_1k\epsilon |\psi|_{C^1(M)} + C_2 |\psi|_{C^2(M)}).
\]

Notice that $|\tilde{B}_{2R_1} (\gamma^\circ_{\tilde{I}})| \leq |B^\circ_{\frac{c}{\epsilon k}}(\Gamma)| \leq \frac{C|M|}{k\epsilon^2}$. Hence

\begin{equation}
\label{eg4}
\int_{\tilde{B}_{R_1} (\gamma^\circ_{\tilde{I}})} \tilde{A} d\psi \leq \frac{C_1}{\epsilon} |\psi|_{C^1(M)} + \frac{C_2}{\epsilon^2 k} |\psi|_{C^2(M)}.
\end{equation}

Estimates (\ref{eg1}), (\ref{eg2}), (\ref{eg3}) and (\ref{eg4}) together imply the estimate (\ref{eg5}). The proof of (\ref{eg6}) is similar and is slightly simpler, since the last step using (2) of proposition \ref{ea} is no longer needed.
\hfill$\Box$\\
\begin{theorem}
\label{ec}
For $\epsilon = \epsilon_k \rightarrow 0$, $\displaystyle \lim_{k \rightarrow +\infty} \frac{1}{k} [\tilde{X}] = \lim_{k \rightarrow +\infty} \frac{1}{k} [\hat{X}] = \omega$ as currents. More precisely, for any smooth $(2n-2)$-form $\psi$ on $M$,
\[
\langle [\tilde{X}], \psi \rangle - k\langle \omega, \psi \rangle \leq C_1k\epsilon_k |\psi|_{C^1(M)} + \frac{C_2}{\epsilon_k^2 k} |\psi|_{C^2(M)},
\]
\[
\langle [\hat{X}], \psi \rangle - k\langle \omega, \psi \rangle \leq Ck\epsilon_k |\psi|_{C^1(M)}.
\]
\end{theorem}
{\bf Proof:} The only difference from theorem \ref{eg} is that for $p \in U_{i,\epsilon}$, $|A_i(p)|_g = O(k|p-p_i|) = O(k\epsilon_k)$ and\\

\hspace{.7in}$\displaystyle|\tilde{A}(p)| \leq C \sum_{j\in I} k|p_i-p_j| e^{-k|p_i-p_j|{\rm Dist}(p,\Gamma_{ij})/2} \leq Ck\epsilon_k.$
\hfill$\Box$\\\\

{\bf Remark:} From theorem \ref{eg}, one can see that our section of $L^k$ is quite different from Donaldson's section constructed in \cite{D1}. In a sense, Donaldson's section is the most generic that use as many peak sections as possible, and our section is less generic that use as few peak sections as possible. Theorem \ref{ec} indicates that our section will behave more like Donaldson's section when $\epsilon_k$ approaches $k^{-1/2}$ scale.\\

\ifx\undefined\bysame
\newcommand{\bysame}{\leavevmode\hbox to3em{\hrulefill}\,}
\fi

\end{document}